\documentclass[letterpaper,11pt]{article}
\usepackage{amsmath}
\usepackage{amsthm}
\usepackage{amsfonts}
\usepackage{amssymb}
\usepackage{latexsym}
\usepackage{amsmath, amsthm, amssymb}
\usepackage{epsfig,graphicx}
\usepackage{graphics}
\usepackage{psfrag}
\usepackage{subfigure}

\usepackage{color}
\usepackage{graphicx}
\usepackage{amscd}
\usepackage{epsfig}

\newtheorem{thm}{Theorem}
\newtheorem{cor}[thm]{Corollary}

\newtheorem{lem}[thm]{Lemma}
\newtheorem{defn}[thm]{Definition}
\newtheorem{prop}[thm]{Proposition}
\newtheorem{clm}[thm]{Proposition}
\newtheorem{rem}[thm]{Remark}
 \numberwithin{equation}{section}
\newenvironment{Proof}{\noindent\bf{Proof.}\rm}{\hfill$\blacksquare$\bigskip}

\newcommand{\whp}{whp}

\newcommand{\PlantedDist}{{\cal{P}}^{{\rm plant}}_{n,p}}

\newcommand{\F}{{\cal{F}}}
\newcommand{\Core}{{\cal{H}}}
\newcommand{\one}{\textbf{1}}
\newcommand{\zero}{\textbf{0}}
\newcommand{\WP}{\textsf{WP}}

\newcommand{\E}{\mathbb{E}}

\definecolor{Red}{rgb}{1,0,0}
\definecolor{Blue}{rgb}{0,0,1}
\definecolor{Olive}{rgb}{0.41,0.55,0.13}
\definecolor{Green}{rgb}{0,1,0}
\definecolor{MGreen}{rgb}{0,0.8,0}
\definecolor{DGreen}{rgb}{0,0.55,0}
\definecolor{Yellow}{rgb}{1,1,0}
\definecolor{Cyan}{rgb}{0,1,1}
\definecolor{Magenta}{rgb}{1,0,1}
\definecolor{Orange}{rgb}{1,.5,0}
\definecolor{Violet}{rgb}{.5,0,.5}
\definecolor{Purple}{rgb}{.75,0,.25}
\definecolor{Brown}{rgb}{.75,.5,.25}
\definecolor{Grey}{rgb}{.5,.5,.5}
\definecolor{Black}{rgb}{0,0,0}

\title{ {Complete convergence of message passing algorithms for
some satisfiability problems } \ } \author{ Uriel
Feige\thanks{The Weizmann Institute. {\tt
uriel.feige@weizmann.ac.il}.} \and Elchanan Mossel\thanks{U.C.
Berkeley and Weizmann Institute. E-mail: {\tt
        mossel@stat.berkeley.edu}. Supported by a Sloan fellowship
  in Mathematics, by NSF Career award DMS-0548249 and Israeli Science Foundation grant 1300/08}  \and Dan Vilenchik\thanks{UCLA. E-mail: {\tt
        vilenchik@math.ucla.edu}.}}

\begin{document}
\maketitle
\begin{abstract}
In this paper we analyze
the performance of \textsf{Warning Propagation}, a popular message
passing algorithm. We show that for 3CNF formulas drawn from a certain
distribution over random satisfiable 3CNF formulas, commonly
referred to as the planted-assignment distribution, running
\textsf{Warning Propagation} in the standard way (run message
passing until convergence, simplify the formula according to the
resulting assignment, and satisfy the remaining subformula, if
necessary, using a simple ``off the shelf" heuristic) results in a satisfying assignment when
the clause-variable ratio is a sufficiently large constant.

\end{abstract}

\newpage

\section{Introduction} A CNF formula over the variables
$x_{1},x_{2},...,x_{n}$ is a conjunction of clauses
$C_{1},C_{2},... ,C_{m}$ where each clause is a disjunction of one
or more literals. Each literal is either a variable or its
negation. A formula is said to be in $k$-CNF form if every clause
contains exactly $k$ literals. A CNF formula is satisfiable if
there is a boolean assignment to the variables such that every
clause contains at least one literal which evaluates to true.
3SAT, the language of all satisfiable 3CNF formulas, is well known
to be NP-complete \cite{Cook71}.

The plethora of worst-case NP-hardness results for many
interesting optimization problems motivates the study of
heuristics that give ``useful'' answers for ``typical'' subset of
the problem instances. In this paper we seek to evaluate those two measures rigorously and for that we shall use random models
and average case analysis.

In this paper we study random satisfiable 3CNF formulas with an arbitrary density. For this we use the {\em planted distribution}
\cite{flaxman,AlonKrivSudCliqe,OnTheGreedy,ExpectedPoly3SAT,AlonKahale97,ChenFrieze} denoted throughout
by $\PlantedDist$. A random 3CNF in this distribution is obtained
by first picking an assignment $\varphi$ to the variables, and
then including every clause satisfied by $\varphi$ with
probability $p=p(n)$, thus guaranteeing that the resulting
instance is satisfiable.


We briefly note that there exists another model of
random 3CNF formulas which consists of $m$ clauses chosen uniformly at random from the set of all
$8\binom{n}{3}$ possible ones ($m$ is a parameter of the distribution, and $m/n$ is referred to as the clause-variable ratio, or the
density, of the random instance). This distribution shows a sharp threshold with respect to satisfiability \cite{Friedgut}.
Specifically, a random 3CNF with clause-variable ratio below the threshold is satisfiable $\whp$
(with high probability, meaning with probability tending to 1 as $n$
goes to infinity) and one with ratio above the threshold is
unsatisfiable $\whp$. Experimental results predict the threshold to be around 4.2 \cite{ExperimentalUpperBoundOnThresh}.
\\\\
To describe our main result, we formally define the \textsf{Warning Propagation (WP)} algorithm.

\subsection{Warning Propagation}\label{WPSubs}
$\WP$ is a simple iterative message passing algorithm similar to
\textsf{Belief Propagation} \cite{Pearl82} and
\textsf{Survey Propagation} \cite{SurveyPropagation}.
Messages in
the $\WP$ algorithm can be interpreted as "warnings", telling
a clause the values that variables will have if the clause "keeps
quite" and does not announce its wishes, and telling a variable
which clauses will not be satisfied if the variable does not commit
to satisfying them. We now present the algorithm in a formal way.

Let $\F$ be a CNF formula. For a variable $x$, let $N^+(x)$ be the
set of clauses in $\F$ in which $x$ appears positively (namely, as
the literal $x$), and $N^{-}(x)$ be the set of clauses in which $x$
appears negatively. For a clause $C$, let $N^+(C)$ be the set of
variables that appear positively in $C$, and respectively $N^-(C)$
for negative ones.

There are two types of messages involved in the
\textsf{WP} algorithm.  Messages sent from a variable
$x_i$ to a clause $C_j$ in which it appears, and vice a versa.
Consider a clause $C_j$ that contains a variable $x_i$. A message from $x_i$ to $C$ is denoted by  $x_i \rightarrow C_j$,  and it has an integer value. A positive value indicates that $x_i$ is tentatively set to true, and a negative value indicates that $x_i$ is tentatively set to false. A message from $C_j$ to $x_i$ is denoted by $C_j \rightarrow x_i$, and it has a Boolean value. A value of 1 indicates that $C_j$ tentatively wants $x_i$ to satisfy it, and a 0 value indicates that $C_j$ is tentatively indifferent to the value of $x_i$. We now present the update rules for these messages.
$$ x_i \rightarrow C_j = \left(\sum_{C_k \in N^+(x_i),k\neq
j} C_k \rightarrow x_i \right) - \left(\sum_{C_k\in N^-(x_i),k\neq j}C_k
\rightarrow x_i\right).$$ If $x_i$ appears only in $C_j$ then we set the
message to 0.
$$ C_j \rightarrow x_i = \left(\prod_{x_k \in N^+(C_j),k\neq
i} I_{<0}(x_k \rightarrow C_j) \right) \cdot \left( \prod_{x_k \in N^-(C_j),k\neq i}
I_{>0}(x_k \rightarrow C_j) \right),$$ where $I_{< 0}(b)$ is an indicator
function which is `1' iff $b<0$ (respectively $I_{>0}$). If $C_j$
contains only $x_i$ (which cannot be the case in 3CNF formulas)
then the message is set to 1.
Lastly, we define the current assignment
of a variable $x_i$ to be
$$ B_i = \left(\sum_{C_j \in N^+(x_i)} C_j \rightarrow x_i \right)- \left(\sum_{C_j\in
N^-(x_i)}C_j \rightarrow x_i\right).$$ If $B_i>0$ then $x$ is assigned
TRUE, if $B_i<0$ then $x_i$ is assigned FALSE, otherwise $x_i$ is
UNASSIGNED. Assume some order on the clause-variable messages
(e.g. the lexicographical order on pairs of the form $(j,i)$
representing the message $C_{j} \rightarrow x_i$). Given a vector
$\alpha \in \{0,1\}^{3m}$ in which every entry is the value of the
corresponding $C_{j}\rightarrow x_i$ message, a partial assignment
$\psi\in \{TRUE,FALSE,UNASSIGNED\}^n$ can be generated according
to the corresponding $B_i$ values (as previously explained).

Given a 3CNF formula $\F$ on $n$ variables and $m$ clauses, the
\textbf{factor graph} of $\F$, denoted by
$FG(\F)$, is the following graph representation of $\F$. The factor
graph is a bipartite graph, $FG(\F)=(V_1\cup V_2,E)$ where
$V_{1}=\{x_{1},x_2,...,x_n\}$ (the set of variables) and
$V_2=\{C_1,C_2,...,C_m\}$ (the set of clauses). $(x_i,C_j)\in E$ iff
$x_i$ appears in $C_j$. For a 3CNF $\F$ with $m$ clauses it holds
that $\#E=3m$, because every clause contains exactly 3 different
variables. (Here and elsewhere, $\#A$ denotes the
cardinality of a set $A$. The notation $|a|$ will denote the absolute value of a real number
$a$.)

It would be convenient to think of the messages in terms of the
corresponding factor graph. Every undirected edge $(x_i,C_j)$ of the
factor graph is replaced with 2 anti-parallel directed edges,
$(x_i\rightarrow C_j)$ associated with the message $x_i\rightarrow
C_j$ and respectively the edge $(C_j\rightarrow x_i)$.

Let us now formally describe the algorithm.
\\\\
$\verb"Warning Propagation(3CNF formula"$ $\F):$ \\
$\verb"1. construct the corresponding factor graph"$ $FG(\F).$\\
$\verb"2. randomly initialize the clause-variable messages to 0 or 1."$\\
$\verb"3. repeat until no clause-variable message changed from the"$\\
$\verb"   previous iteration:"$\\
$\verb"   3.a randomly order the edges of" $ $FG(\F).$\\
$\verb"   3.b update all clause-variable messages"$
$C_j\rightarrow
x_i$ $\verb"according"$\\
$\verb"       to the random edge order."$\\
$\verb"4. compute a partial assignment"$ $\psi$ $\verb"according to the"$ $B_i$ $\verb"messages".$\\
$\verb"5. return "$ $\psi$.\\


The variable-clause message updates are implicit in line
3.b: when evaluating the clause-variable message
along the edge $C \rightarrow x$, $C=(x\vee y\vee z)$, the
variable-clause messages concerning this calculation ($z,y
\rightarrow C$) are evaluated on-the-fly using the last updated
values $C_i \rightarrow y$, $C_j \rightarrow z$ (allowing feedback
from the same iteration). We allow the algorithm not to terminate (the clause-variable messages may keep changing every
iteration). If the algorithm does return an assignment $\psi$ then
we say that it converged. In practice it is common to limit in
advance the number of iterations, and if the algorithm does not
converge by then, return a failure.

\subsection{Main Results}\label{ResutlsSubs}
Our contribution is analyzing the performance of \textsf{Warning
Propagation} (\textsf{WP} for brevity), a popular message passing
algorithm, when applied to satisfiable formulas drawn from a certain
random distribution over satisfiable 3CNF formulas, commonly called
the planted distribution. We show that the standard way of running
message passing algorithms -- run message passing until convergence,
simplify the formula according to the resulting assignment, and
satisfy the remaining subformula, if possible, using a simple ``off
the shelf" heuristic -- works for planted random satisfiable
formulas with a sufficiently large constant clause-variable ratio.
As such, our result is the first to rigorously prove the
effectiveness of a message passing algorithm for the solution of a
non-trivial random SAT distribution.
We note that a recent work \cite{AminAchi} demonstrated the usefulness of analytical tools developed for the planted distribution
for ``hard" uniform instances.
To formally state our result we require a few additional definitions.

Given a 3CNF $\F$, \textbf{simplify} $\F$ according to $\psi$,
when $\psi$ is a partial assignment, means: in every clause
substitute every assigned variable with the value given to it by
$\psi$. If a clause contains a literal which evaluates to
true, remove the clause. From the remaining clauses, remove all literals which
evaluate to false. The resulting instance is not necessarily in
3CNF form, as clauses may have any number of literals between~0 and~3.
Denote by $\F|_{\psi}$ the 3CNF $\F$ simplified
according to $\psi$. Note that $\F|_{\psi}$ may contain empty clauses, in which case it is not satisfiable.
For a set of variables $A\subseteq V$, denote
by $\F[A]$ the set of clauses in which all variables belong to
$A$.

We call a 3CNF formula \textbf{simple}, if it can be satisfied using simple
well-known heuristics (examples include very sparse random 3CNF
formulas which are solvable $\whp$ using the pure-literal
heuristic \cite{BorderFriezeUpfal}, formulas with small weight
terminators -- to use the terminology of \cite{BenSassonAlke} --
solvable $\whp$ using RWalkSat, etc). This is a somewhat informal notion, but
it suffices for our needs.

\begin{thm}\label{ConvergenceThmMed} Let $\F$ be a 3CNF formula
randomly sampled according to $\PlantedDist$, where $p \geq
d/n^2$, $d$ a sufficiently large constant. Then the following
holds $\whp$ (the probability taken over the choice of $\F$, and
the random choices in lines 2 and 4 of the \textsf{WP} algorithm).
There exists a satisfying assignment $\varphi^*$ (not necessarily
the planted one) such that:

\begin{enumerate}
  \item  $\textsf{WP}(\F)$ converges after at most $O(\log n)$
  iterations.
  \item  Let $\psi$ be the partial assignment returned by $\textsf{WP}(\F)$, let
  $V_A$ denote the variables assigned to either TRUE or FALSE in
  $\psi$, and $V_U$ the variables left UNASSIGNED. Then for every
  variable $x\in V_A$,
  $\psi(x)=\varphi^*(x)$. Moreover, $\#{V_A} \geq (1-e^{-\Theta(d)})n$.
  \item
  $\F|_{\psi}$ is a simple formula which can be satisfied in time $O(n)$.
\end{enumerate}
\end{thm}

\begin{rem}\rm
Theorem~\ref{ConvergenceThmMed} relates to the planted 3SAT model,
but \cite{UniformSAT} implies that it also true to the
uniform random 3SAT distribution, in which a \emph{satisfiable} formula with $m$ clauses
is chosen uniformly at random among all such formulas.
\end{rem}

\begin{prop}\label{ConvergencePropDense} Let $\F$ be a 3CNF formula
randomly sampled according to $\PlantedDist$, where $p\geq c\log n
/n^2$, with $c$ a sufficiently large constant, and let $\varphi$ be
its planted assignment. Then $\whp$ after at most 2 iterations,
$\textsf{WP}(\F)$ converges, and the returned $\psi$ equals
$\varphi$.
\end{prop}

It is worth noting that formulas in $\PlantedDist$, with $n^2p$
some large constant, are \emph{not} known to be simple (in the sense that
we alluded to above). For example, it is shown in
\cite{BenSassonAlke} that RWalkSat is very unlikely to hit a
satisfying assignment in polynomial time when running on a random
$\PlantedDist$ instance in the setting of
Theorem \ref{ConvergenceThmMed}. Nevertheless, planted 3CNF formulas with sufficiently large
(constant) density were shown to be solvable $\whp$ in $\cite{flaxman}$ using a spectral algorithm.
Though in our analysis we use similar techniques to
\cite{flaxman} (which relies on \cite{AlonKahale97}), our result is conceptually different in the
following sense. In \cite{AlonKahale97,ChenFrieze,flaxman}, the starting point is the planted distribution, and then one designs an
algorithm that works well under this distribution. The algorithm may
be designed in such a way that makes its analysis easier. In
contrast, our starting point is a given algorithm, $\WP$, and then we ask for which input distributions it
works well. We cannot change the algorithm in ways that would
simplify the analysis. Another difference between our work and that of
\cite{AlonKahale97,ChenFrieze,flaxman} is that unlike the algorithms
analyzed in those other papers, \textsf{WP} is a randomized
algorithm, a fact which makes its analysis more difficult. We could
have simplified our analysis had we changed \textsf{WP} to be
deterministic (for example, by initializing all clause-variable
messages to~1 in step~2 of the algorithm), but there are good
reasons why \textsf{WP} is randomized. For example, it can be shown
that (the randomized version) \textsf{WP} converges with
probability~1 on 2CNF formulas that form one cycle of implications,
but might not converge if step~4 does not introduce fresh randomness
in every iteration of the algorithm.

\medskip

The planted 3SAT model is also similar to LDPC codes in many ways. Both constructions are
based on random factor graphs. In codes, the received corrupted
codeword provides noisy information on a single bit or on the
parity of a small number of bits of the original codeword. In
$\PlantedDist$, $\varphi$ being the planted assignment, the
clauses containing a variable $x_i$ provide noisy information on
the polarity of $\varphi(x_i)$.
Comparing our results with the coding setting, the effectiveness of
message passing algorithms for amplifying local information in order
to decode codes close to channel capacity was established
in a number of papers, e.g.~\cite{LMSS98,RSU01}. Our results are
similar in flavor, however the combinatorial analysis provided here
allows to recover an assignment satisfying \emph{all} clauses,
whereas in the random LDPC codes setting, message passing allows to
recover only $1-o(1)$ fraction of the codeword correctly. In
\cite{LMSS01} it is shown that for the erasure channel, all bits may
be recovered correctly using a message passing algorithm, however in
this case the LDPC code is designed so that message passing works
for it.

\medskip

It is natural to ask whether our analysis can be extended to
show that \textsf{Belief Propagation} (\textsf{BP}) finds a
satisfying assignment to $\PlantedDist$ in the setting of Theorem
\ref{ConvergenceThmMed}. Experimental results predict the answer
to be positive. However our analysis of \textsf{WP} does not
extend as is to \textsf{BP}. In \textsf{WP}, all warnings received
by a variable (or by a clause) have equal weight, but in
\textsf{BP} this need not be the case (there is a probability
level associated with each warning). In particular, this may lead
to the case where a small number of ``very" wrongly assigned variables
``pollute" the entire formula, a possibility that our analysis
managed to exclude for the \textsf{WP} algorithm.

\subsection{Paper's organization} The remainder of the paper is structured as follows. Section
\ref{sec:intuition} provides an overview that may help the reader
follow the more technical parts of the proofs. In Section
\ref{PropOfRandomInstSection} we discuss some properties that a
typical instance in $\PlantedDist$ possesses. Using these
properties, we prove in Section \ref{ResultsProofSection} Theorem
\ref{ConvergenceThmMed} and Proposition
\ref{ConvergencePropDense}.

\section{Proof Outline} \label{sec:intuition}

Let us first consider some possible fixed points of the Warning
Propagation (WP) algorithm. The {\em trivial} fixed point is the
one in which all messages are \zero. One may verify that this is
the unique fixed point in some cases when the underlying 3CNF
formula is very easy to satisfy, such as when all variables appear
only positively, or when every clause contains at least two
variables that do not appear in any other clause. A {\em local
maximum} fixed point is one that corresponds to a strict local
maximum of the underlying MAX-3SAT instance, namely to an
assignment $\tau$ to the variables in which flipping the truth
assignment of any single variable causes the number of satisfied
clauses to strictly decrease. The reader may verify that if every
clause $C$ sends a \one\ message to a variable if no other
variable satisfies $C$ under $\tau$, and a \zero\ message
otherwise, then this is indeed a fixed point of the WP algorithm.
Needless to say, the WP algorithm may have other fixed points, and
might not converge to a fixed point at all.

Recall the definition of $\PlantedDist$. First a truth assignment
$\varphi$ to the variables $V=\{x_{1},x_{2},...,x_{n}\}$ is picked
uniformly at random. Next, every clause satisfied by $\varphi$ is
included in the formula with probability $p$ (in our case $p\geq
d/n^2$, $d$ a sufficiently large constant). There are
$(2^3-1)\cdot\binom{n}{3}$ clauses satisfied by $\varphi$, hence
the expected size of $\F$ is $p\cdot7\cdot\binom{n}{3} = 7 d
n/6+o(n)$ (when $d$ is constant, then this is linear in $n$, and
therefore such instances are sometimes referred to as
\emph{sparse} 3CNF formulas).  To simplify the presentation, we
assume w.l.o.g. (due to symmetry) that the planted assignment
$\varphi$ is the all-one vector.

To aid intuition, we list some (incorrect) assumptions and analyze
the performance of WP on a $\PlantedDist$ instance under these
assumptions.

\begin{enumerate}

\item In expectation, a variable appears in
$4\binom{n}{2}p=2d+o(1)$ clauses positively, and in $3d/2+o(1)$
clauses negatively. Our first assumption is that for every
variable, its number of positive and negative appearances is equal
to these expectations.

\item We say that a variable {\em supports} a clause with respect
to the planted assignment (which was assumed without loss of
generality to be the all \one\ assignment) if it appears
positively in the clause, and the other variables in the clause
appear negatively. Hence the variable is the only one to satisfy
the clause under the planted assignment. For every variable in
expectation there are roughly $d/2$ clauses that it supports. Our
second assumption is that for every variable, the number of
clauses that it supports is equal to this expectation.

\item Recall that in the initialization of the \textsf{WP}
algorithm, every clause-variable message $C\to x$ is 1 w.p.
$\frac{1}{2}$, and 0 otherwise. Our third assumption is that with
respect to every variable, half the messages that it receives from
clauses in which it is positive are initialized to~1, and half the
messages that it receives from clauses in which it is negative are
initialized to~1.

\item Recall that in step 3b of \textsf{WP}, clause-variable
messages are updated in a random order. Our fourth assumption is
that in each iteration of step~3, the updates are based on the
values of the other messages from the previous iteration, rather
than on the last updated values of the messages (that may
correspond either to the previous iteration or the current
iteration, depending on the order in which clause-variable
messages are visited). Put differently, we assume that in step 3b
all clause-variable messages are evaluated in \emph{parallel}.

\end{enumerate}

Observe that under the first two assumptions, the planted
assignment is a local maximum of the underlying MAX-3SAT instance.
We show that under the third and fourth assumption, \textsf{WP}
converges to the corresponding local maximum fixed point in two
iterations: Based on the initial messages as in our third
assumption, the messages that variables send to clauses are all
roughly $(2d-3d/2)/2=d/4$. Following the initialization, in the
first iteration of step~3 every clause $C$ that $x$ supports will
send $x$ the message 1, and all other messages will be 0. Here we
used our fourth assumption. (Without our fourth assumption,
\textsf{WP} may run into trouble as follows. The random ordering
of the edges in step 3 may place for some variable $x$ all
messages from clauses in which it appears positively before those
messages from clauses in which it appears negatively. During the
iteration, some of the messages from the positive clauses may
change from 1 to 0. Without our fourth assumption, this may at
some point cause $x$ to signal to some clauses a negative rather
than positive value.) The set of clause-variable messages as above
will become a fixed point and repeat itself in the second
iteration of step 3. (For the second iteration, the fourth
assumption is no longer needed.) Hence the algorithm will
terminate after the second iteration.

Unfortunately, none of the four assumptions that we made are
correct. Let us first see to what extent they are violated in the
context of Proposition~\ref{ConvergencePropDense}, namely, when
$d$ is very large, significantly above $\log n$. Standard
concentration results for independent random variables then imply
that the first, second and third assumptions simultaneously hold
for all variables, up to small error terms that do not effect the
analysis. Our fourth assumption is of course never true, simply
because we defined \textsf{WP} differently. This complicates the
analysis to some extent and makes the outcome depend on the order
chosen in the first iteration of step 3a of the algorithm.
However, it can be shown that for most such orders, the algorithm
indeed converges to the fixed point that corresponds to the
planted assignment.

The more difficult part of our work is the case when $d$ is
a large constant. In this case, already our first
two assumptions are incorrect. Random fluctuations with respect to
expected values will $\whp$ cause a linear fraction of the
variables to appear negatively more often than positively, or not
to support any clause (with respect to the planted assignment). In
particular, the planted assignment would no longer be a local
maximum with respect to the underlying MAX-3SAT instance.
Nevertheless, as is known from previous work~\cite{flaxman}, a
large fraction of the variables will behave sufficiently close to
expectation so that the planted assignment is a local maximum with
respect to these variables. Slightly abusing notation, these set
of variables are often called the {\em core} of the 3CNF formula.
Our proof plan is to show that $\WP$ does converge, and that the
partial assignment in step 4 assigns all core variables their
correct planted value. Moreover, for non-core variables, we wish
to show that the partial assignment does not make any
unrecoverable error -- whatever value it assigns to some of them,
it is always possible to assign values to those variables that are
left unassigned by the partial assignment so that the input
formula is satisfied. The reason why we can expect such a proof
plan to succeed is that it is known to work if one obtains an
initial partial assignment by means other than $\WP$, as was already
done in~\cite{flaxman,TechReport}.

Let us turn now to our third assumption. It too is violated for a
linear fraction of the variables, but is nearly satisfied for most
variables. This fact marks one point of departure for our work
compared to previous work~\cite{flaxman,TechReport}. Our
definition of the core variables will no longer depend only on the
input formula, but also on the random choice of initialization
messages. This adds some technical complexity to our proofs.

The violation of the fourth assumption is perhaps the technical
part in which our work is most interesting. It relates to the
analysis of $\WP$ on factor graphs that contain cycles, which is
often a stumbling point when one analyzes message passing
algorithms. Recall that when $d$ is very large
(Proposition~\ref{ConvergencePropDense}), making the fourth
assumption simplifies the proof of convergence of WP. Hence
removing this assumption in that case becomes a nuisance. On the
other hand, when $d$ is smaller (as in
Theorem~\ref{ConvergenceThmMed}), removing this assumption becomes
a necessity. This will become apparent when we analyze convergence
of WP on what we call {\em free cycles}. If messages in step 3b of
\textsf{WP} are updated based on the value of other messages in
the {\em previous} iteration (as in our fourth assumption), then
the random choice of order in step 3a of \textsf{WP} does not
matter, and one can design examples in which the messages in a
free cycle never converge. In contrast, if messages in step 3b of
\textsf{WP} are updated based on the latest value of other
messages (either from the previous iteration or from the current
iteration, whichever one is applicable), free cycles converge with
probability 1 (as we shall later show).

To complete the proof plan, we still need to show that simplifying
the input formula according to the partial assignment returned by
\textsf{WP} results in a formula that is satisfiable, and
moreover, that a satisfying assignment for this sub-formula can
easily be found. The existential part (the sub-formula being
satisfiable) will follow from a careful analysis of the partial
assignment returned by \textsf{WP}. The algorithmic part (easily
finding an assignment that satisfies the sub-formula) is based on
the same principles used in~\cite{AlonKahale97,flaxman}, showing
that the sub-formula breaks into small connected components.

\section{Properties of a Random $\PlantedDist$
  Instance}\label{PropOfRandomInstSection}

In this section we discuss relevant properties of a random
$\PlantedDist$ instance. This section is rather technical in nature. The proofs
are based on probabilistic arguments that are standard in our context. In the rest of the paper,
for simplicity of presentation, we assume w.l.o.g. that the planted assignment is the all TRUE assignment.

\subsection{Preliminaries}

The following well-known concentration result (see, for example
 \cite[p. 21]{JRL2000}) will be used several times in the
proof. We denote by $B(n,p)$ the binomial random variable with
parameters $n$ and $p$, and expectation $\mu=np$.
\begin{thm}\label{Thm:Chernoff2}(Chernoff's inequality) If $X \sim B(n,p)$ and $t \geq 0$ is some number, then
\[ Pr\big(X \geq \mu + t) \leq e^{-\mu \cdot f(t/\mu)}, \]
\[ Pr\big(X \leq \mu - t) \leq e^{-\mu \cdot f(-t/\mu)},\]

where $f(x) = (1+x)\ln (1+x) - x$.
\end{thm}
\noindent We use the following Azuma-Hoeffding inequality for martingales:
\begin{thm}\label{Thm:RegMartingale} Let $\{X_i\}_{i=0}^{N}$ be a martingale with $|X_k - X_{k-1}| \leq c_k$, then
\[Pr \left[|X_N - X_0| \geq t\right] \leq 2 e^{-t^2 / \sum_{k=1}^N c_k^2}\]
\end{thm}
We shall also use the following version of the martingale argument, which can be found in
\cite[Page 101]{TheProbMethod}. Assume the probability space is generated by a finite set of mutually
independent Yes/No choices, indexed by $i\in I$. Given a random variable $Y$
on this space, let $p_i$ denote the probability that choice $i$ is Yes. Let $c_i$ be such that
changing choice $i$ (keeping all else the same) can change $Y$ by at most $c_i \leq c$.
We call $p_i (1 - p_i )c^2$ the variance of choice $i$. Now consider a solitaire game in which one finds the value of $Y$ by making queries of an always truthful oracle.
The choice of query can depend on previous responses. All possible questioning lines
can be naturally represented in a decision tree form. A ``line of questioning''
is a path from the root to a leaf of this tree, a sequence of questions and responses that
determine $Y$. The total variance of a line of questioning is the sum of the variances
of the queries in it.

\begin{thm}\label{thm_martingale} For every $\varepsilon > 0$ there exists a $\delta > 0$ so that
if for every line of questioning, with parameters
$p_1,p_2,\ldots,p_n$ and $c_1,c_2,\ldots,c_n$, that determines $Y$,
the total variance is at most $\sigma^2$, then
$$Pr[|Y -E[Y ]| > \alpha\sigma] \leq  2e^{-\alpha^2/(2(1+\varepsilon))},$$
for all positive $\alpha$ with $\alpha c < \sigma(1 +
\varepsilon)\delta$, where $c$ is such that $\forall i \phantom{i} c_i \leq c$.
\end{thm}

\subsection{Stable Variables}\label{StableVarsSubs}
\begin{defn} A variable $x$ \textbf{supports} a clause $C$
with respect to a partial assignment $\psi$, if it is the only
variable to satisfy $C$ under $\psi$, and the other two variables
are assigned by $\psi$. \end{defn}
\begin{prop}\label{SupportSuccRate} Let $\F$ be a 3CNF formula
randomly sampled according to $\PlantedDist$, where $p = d/n^2$,
$d \geq d_0$, $d_0$ a sufficiently large constant. Let $F_{SUPP}$ be a random
variable counting the number of variables in $\F$ whose support
w.r.t. $\varphi$ is less than $d/3$. Then $\whp$ $F_{SUPP}\leq
e^{-\Theta(d)}n$.
\end{prop}
\begin{Proof}
Every variable is expected to support
$\frac{d}{n^2}\cdot\binom{n}{2}=\frac{d}{2}+O(\frac{1}{n})$ clauses,
thus using Chernoff's inequality the probability of a single variable supporting less than $d/3$ clauses
is at most $e^{-\Theta(d)}$. Using linearity of expectation, $\mathbb{E}[F_{SUPP}]\leq e^{-\Theta(d)}n$. To prove concentration
around the expected value we use Chernoff's bound once more
as the support of one variable is independent of the others (since
it concerns different clauses which are included independently of
each other). The claim then follows.
\end{Proof}

Following the definitions in Section \ref{WPSubs}, given a 3CNF $\F$ with
a satisfying assignment $\varphi$, and a variable $x$, we let $N^{++}(x)$ be the set of clauses in $\F$
in which $x$ appears positively but doesn't support w.r.t.
$\varphi$. Let $N^{s}(x)$ be the set of clause in $\F$ which $x$
supports w.r.t. $\varphi$. Let $\pi=\pi(\F)$ be some ordering of the
clause-variable message edges in the factor graph of $\F$. For an
index $i$ and a literal $\ell_x$ (by $\ell_x$ we denote a literal
over the variable $x$) let $\pi^{-i}(\ell_x)$ be the set of
clause-variable edges $(C \rightarrow x)$ that appear before index
$i$ in the order $\pi$ and in which $x$ appears in $C$ as $\ell_x$.
For a set of clause-variable edges $\cal{E}$ and a set of clauses
$\cal{C}$ we denote by ${\cal{E}}\cap{\cal{C}}$ the subset of edges
containing a clause from $\cal{C}$ as one endpoint.

\begin{defn}\label{StableDefn} Let $\pi$ be an ordering of the clause-variable messages of a 3CNF formula
$\F$. Let $\varphi$ be a satisfying assignment of $\F$. A variable $x$ is \textbf{$d$-stable} in $\F$ w.r.t. $\pi$ and $\varphi$ if for every location $i$ in $\pi$ that contains
a message $C\rightarrow x$, $C=(\ell_x\vee \ell_y \vee
\ell_z)$, the following holds:
\begin{enumerate} \item  $|\#\pi^{-i}(y)\cap
N^{++}(y)-\#\pi^{-i}(\bar{y})\cap N^{-}(y)|\leq d/30$.
\item $|\#N^{++}(y)-\#N^{-}(y)|\leq d/30$.
\item $\#N^{s}(y)\geq d/3$
\end{enumerate}
and the same holds for $z$.
\end{defn}
When $d$ is clear from context, which will usually be the case, we will suppress the $d$ in the ``$d$-stable".
\begin{prop}\label{StableSuccRate} Let $\F$ be a 3CNF formula
randomly sampled according to $\PlantedDist$, where $p = d/n^2$,
$d \geq d_0$, $d_0$ a sufficiently large constant. Let $\pi$ be a random ordering of
the clause-variable messages, and $F_{UNSTAB}$ be a random variable
counting the number of variables in $\F$ which are \emph{not} stable. Then
$\whp$ $F_{UNSTAB}\leq e^{-\Theta(d)}n$.
\end{prop}
\begin{Proof}
We start by bounding $\mathbb{E}[F_{UNSTAB}]$. Consider a
clause-variable message edge $C \rightarrow x$ in location $i$ in
$\pi$, $C=(\ell_x\vee \ell_y \vee \ell_z)$. Now consider location $j
\leq i$. The probability of an edge $C' \rightarrow \bar{y}$ in
location $j$ is
$\left(3\binom{n}{2}\right)/\left(7\binom{n}{3}\right)=\frac{9}{7n}+O(\frac{1}{n^2})$
which is exactly the probability of an edge $C'' \rightarrow y$,
$C'' \in N^{++}(y)$. This implies
$$E[|\#\pi^{-i}(y)\cap N^{++}(y)-\#\pi^{-i}(\bar{y})\cap
N^{-}(y)|]=0.$$ If however
$$|\#\pi^{-i}(y)\cap N^{++}(y)-\#\pi^{-i}(\bar{y})\cap
N^{-}(y)| > d/30$$ then at least one of the quantities deviates
from its expectation by $d/60$.

Look at $\#\pi^{-i}(y)\cap N^{++}(y)$ -- this is the number of
successes in draws without replacement. It is known that this
quantity is more concentrated than the corresponding quantity if the
draws were made with replacement~\cite{Hoeffding:63}. In particular,
since the expectation of $\#\pi^{-i}(y)\cap N^{++}(y)$ is $O(d)$ it
follows from Chernoff's bound that the probability that it deviates
from its expectation by more than $d/60$ is $e^{-\Theta(d)}$. A
similar statement holds for $\#\pi^{-i}(\bar{y})\cap N^{-}(y)$.
Properties $(b)$ and $(c)$ are bounded similarly using concentration
results.

\noindent The calculations above hold in particular for the first
$5d$ appearances of messages involving $x$. As for message $5d+1$,
the probability of this message causing $x$ to become unstable is
bounded by the event that $x$ appears in more than $5d$ clauses. As
$x$ is expected to appear in $3.5d$ clauses, the latter event
happens w.p. $e^{-\Theta(d)}$ (again using Chernoff's bound). To sum up,
$$Pr[x\text{ is unstable }]\leq 5d\cdot
e^{-\Theta(d)}+e^{-\Theta(d)}=e^{-\Theta(d)}.$$ The bound on
$E[F_{UNSTAB}]$ follows by linearity of expectation.

We are now left with proving that $F_{UNSTAB}$ is concentrated
around its expectation, we do so using a martingale argument.
Define two new random variables, $F_1$ counting the number of
unstable variables $x$ s.t. there exists a clause $C$, containing
$x$, and another variable $y$, s.t. $y$ appears in more than $\log
n$ clauses, and $F_2$ to be the unstable variables s.t. in all
clauses in which they appear, all the other variables appear in at
most $\log n$ clauses. Observe that $F_{UNSTAB}=F_1+F_2$. To bound
$F_1$, observe that if $F_1\geq 1$, then in particular this
implies that there exists a variable which appears in more than
$\log n$ clauses in $\F$. This happens with probability $o(1)$
since every variable is expected to appear only in $O(d)$
clauses (using Chernoff's bound). To bound $F_2$ we use a martingale argument in the
constellation of \cite{TheProbMethod}, page 101. We use the
clause-exposure martingale (the clause-exposure martingale
implicitly includes the random ordering $\pi$, since one can think
of the following way to generate the random instance -- first
randomly shuffle all possible clauses, and then toss the coins).
The exposure of a new clause $C$ can change $F_2$ by at most
$6\log n$ since every variable in $C$ appears in at most $\log n$
clauses, namely with at most $2\log n$ other variables that might
become (un)stable due to the new clause. The martingale's total
variance, in the terminology of Theorem \ref{thm_martingale}, is
$\sigma^2=\Theta(dn\log ^2 n)$. Theorem \ref{thm_martingale}, with
$\alpha=e^{-\Theta(d)}\sqrt{n}/\log n$, and the fact that
$E[F_2]\leq E[F_{UNSTAB}]$, implies concentration around the expectation
of $F_2$.
\end{Proof}

Let $\alpha\in\{0,1\}^{3\#\F}$ be a clause-variable message vector.
For a set of clause-variable message edges $\cal{E}$ let
$\one_\alpha(\cal{E})$ be the set of edges along which the value is
1 according to $\alpha$. For a set of clauses $\cal{C}$,
$\one_\alpha(\cal{C})$ denotes the set of clause-variable message
edges in the factor graph of $\F$ containing a clause from $\cal{C}$
as one endpoint and along which the value is 1 in $\alpha$.

\begin{defn}\label{ViolatedDefn} Let $\pi$ be an ordering of the clause-variable messages of a 3CNF formula
$\F$, and $\alpha$ a clause-variable message vector. Let $\varphi$ be a satisfying assignment of $\F$.
We say that a variable $x$ is \textbf{$d$-violated} by $\alpha$ in $\pi$ if there exists a
message $C \rightarrow x$, $C=(\ell_x\vee \ell_y\vee \ell_z)$, in
place $i$ in $\pi$ s.t. one of the following holds:
\begin{enumerate}
\item  $|\#\one_\alpha(\pi^{-i}(y)\cap
N^{++}(y))-\#\one_\alpha(\pi^{-i}(\bar{y})\cap N^{-}(y))|> d/30$
\item $|\#\one_\alpha(N^{++}(y))-\#\one_\alpha(N^{-}(y))|>
d/30$
\item $\#\one_\alpha(N^s(y))< d/7$.
\end{enumerate}
Or one of the above holds for $z$.
\end{defn}

We say that a variable is \emph{$r$-bounded} in $\F$ if it appears in at most $r$ clauses.
We say that a variable $x$ has an \emph{$r$-bounded neighborhood} in $\F$ if every clause $C=(\ell_x\vee \ell_y \vee
\ell_z)$ in $\F$ that contains $x$ is such that $y$ and $z$ are $r$-bounded ($x$ itself is not limited)

\begin{lem}\label{GeneralViolatedSuccRate} 
Let F be an arbitrary satisfiable 3CNF formula, 
let $\psi$ be an arbitrary satisfying assignment for $F$,  
and let $\pi$ be an arbitrary ordering of clause-variable messages. 
For a given value $d$, 
let $X$ be the set of variables that have a $20d$-bounded neighborhood in $F$ and are $d$-stable with respect to $\psi$ and $\pi$.
Let $\alpha$ be a random clause-variable message
vector, and $F_{VIO}$ a random variable counting the number of
$d$-violated variables in $X$. Then $\whp$ $F_{VIO}\leq e^{-\Theta(d)}n$.
\end{lem}
\begin{Proof}
First observe that the probability in the statement is taken only over the coin tosses in the choice of $\alpha$, as all other parameters are fixed.
As in the proof of Proposition
\ref{StableSuccRate}, we first bound $\mathbb{E}[F_{VIO}]$, and then
prove concentration using a martingale argument. Fix $x \in X$ and a clause-variable message $C
\rightarrow x$ at location $i$ in $\pi$, $C=(\ell_x\vee \ell_y \vee
\ell_z)$. Let $A^+=\pi^{-i}(y)\cap
N^{++}(y)$ (the set of messages preceding location $i$ where $y$
appears positively but doesn't support the clause) and
$A^-=\pi^{-i}(y)\cap N^{-}(y)$ (the set of messages preceding
location $i$ where $y$ appears negatively). Let us assume w.l.o.g
that $\#A^+ \geq \#A^-$. Stability implies
$$\#A^+-\#A^-\leq d/30.$$
Since $\alpha$ is a random assignment to the clause-variable
messages,
$$\E[\#\one_\alpha(A^+)]-\E[\#\one_\alpha(A^-)]\leq
d/60.$$ If however
\begin{equation}\label{equa}
|\#\one_\alpha(A^+)-\#\one_\alpha(A^-)|> d/30,
\end{equation}
then at least one of $\#\one_\alpha(A^+),\#\one_\alpha(A^-)$
deviated from its expectation by at least $(d/30-d/60)/2=d/120$.
Both quantities are binomially distributed with expectation
$O(d)$ ($x$ has a $20d$-bounded neighborhood, and therefore $y$ appears in at most $20d$ clauses), and therefore the probability of the latter happening is $e^{-\Theta(d)}$ (using standard concentration results). Properties $(b)$ and
$(c)$ are bounded similarly, and the same argument takes care of $z$. Using the
union bound and the linearity of expectation, one obtains that $\E[F_{VIO}]
\leq e^{-\Theta(d)}\#X \leq e^{-\Theta(d)}n$. Let us use Theorem \ref{Thm:RegMartingale} to obtain concentration around $\E[F_{VIO}]$.
Let $Y$ be the set of variables that share some clause with a variable in $X$. Expose the value of messages $C \to y$ for $y \in Y$. Since all $y \in Y$ are $20d$-bounded, the length of the martingale, $N = O(dn)$. The boundedness condition also gives that the martingale difference is $O(d)$. Taking $t=e^{-\Theta(d)}n$, and plugging all the parameters in Theorem \ref{Thm:RegMartingale}, the result follows.

\end{Proof}

\begin{prop}\label{ViolatedSuccRate} Let $\F$ be a 3CNF formula
randomly sampled according to $\PlantedDist$, where $p = d/n^2$,
$d \geq d_0$, $d_0$ is a sufficiently large constant. Let $\pi$ be a random ordering of
the clause-variable messages and $\alpha$ a random clause-variable message vector.
Then $\whp$ the number of stable variables which are violated in $\alpha$ is at most $e^{-\Theta(d)}n$.
\end{prop}
\begin{Proof}
Let $X'$ be the set of stable variables in $\F$ w.r.t.~$\pi$, and $X \subseteq X'$ be the set of variables with a $20d$-bounded
neighborhood in $\F$. Lemma \ref{GeneralViolatedSuccRate} guarantees that the number of violated variables in $X$ is
at most $e^{-\Theta(d)}n$.
It suffices to show that $\#(X' \setminus X) \leq e^{-\Theta(d)}n$. Let $Z_t$ be the set of variables that appear in $t$ clauses in $\F$. The number of clauses in which a variable appears is binomially distributed with expected value $\mu = 7\binom{n}{2}d/n^2 \leq 7d/2$.
Using Theorem \ref{Thm:Chernoff2}, it holds that
$$Pr[x \in Z_t] \leq e^{-\mu \cdot f((t-\mu)/\mu)},$$
where $f(x) = (1+x)\ln (1+x) - x$. For $t \geq 20d$, $f((t-\mu)/\mu) \geq t/(2\mu)$, and therefore
$$Pr[x \in Z_t \, , t \geq 20d] \leq e^{-t/2}.$$
Using the linearity of expectation, and standard concentration results, it holds that $\whp$ for every $t \geq 20d$, $\# Z_t \leq e^{-t/3}n$.
The number of variables in $X' \setminus X$ is then $whp$ at most
$$\sum_{t \geq 20d} 2t \cdot \#Z_t \leq \sum_{t \geq 20d} 2t \cdot e^{-t/3}n \leq e^{-\Theta(d)}n.$$
(Every variable in $Z_t$ ``spoils" at most two other variables in every clause in which it appears, which possibly end up in $X' \setminus X$).
\end{Proof}

\subsection{Dense Subformulas}\label{DesnseSubSubs}
The next property we discuss is analogous to a property proved in
\cite{AlonKahale97} for random graphs. Loosely speaking,
a random graph typically doesn't
contain a small induced subgraph with a large average degree. A
similar proposition for 3SAT can also be found in \cite{flaxman}.
\begin{prop}\label{NoDenseSubformulas} Let $c > 1$ be an arbitrary constant. Let $\F$ be a 3CNF formula
randomly sampled according to $\PlantedDist$, where $p = d/n^2$,
$d \geq d_0$, $d_0$ a sufficiently large constant. Then $\whp$ there exists \emph{no} subset
of variables $U$, s.t. $\#U \leq e^{-\Theta(d)} n$ and there are
at least $c\#U$ clauses in $\F$ containing two variables from $U$.
\end{prop}
\begin{Proof}
For a fixed set $U$ of variables, $\#U=k$, the number of clauses
containing
two variables from $U$ is $$\binom{k}{2}(n-2)2^3 \leq 4k^2n.$$\\
Each of these clauses is included independently w.p.
$\frac{d}{n^2}$. Thus, the probability that $ck$ of them are
included is at most
$$\binom{4k^2n}{ck}\left(\frac{d}{n^2}\right)^{ck}\leq
\left(\frac{4k^2ne}{ck}\cdot \frac{d}{n^2}\right)^{ck} \leq
\left(\frac{12
 kd }{cn}\right)^{ck}.$$
Using the union bound, the probability there exists a ``dense" set
$U$ is at most
$$\sum_{k=2}^{e^{-\Theta(d)} n}\binom{n}{k}\left(\frac{12
 kd }{cn}\right)^{ck}=O(d^{2c}/n^{2c-2}).$$
The last equality is obtained using standard calculations, and the standard estimate on the binomial coefficient:
$\binom{n}{k} \leq (en/k)^k$.
\end{Proof}

\subsection{The Core Variables}\label{CoreSubsection}
We describe a subset of the variables, denoted throughout by
$\Core$ and referred to as the \emph{core variables}, which plays
a crucial role in the analysis. The notion of a stable variable is
not enough to ensure that the algorithm will set a stable variable
according to the planted assignment, as it may happen that a
stable variable $x$ appears in many of its clauses with  unstable
variables. Thus, $x$ can be biased in the wrong direction (by
wrong we mean disagreeing with the planted assignment). However,
if most of the clauses in which $x$ appears contain only stable
variables, then this is already a sufficient condition to ensure
that $x$ will be set correctly by the algorithm. The set $\Core$
captures the notion of such variables. There are several ways to
define a set of variables with these desired properties, we
present one of them, and give a constructive way of obtaining it
(though it has no algorithmic implications, at least not in our
context).

\medskip

\noindent Formally, $\Core=\Core(\F,\varphi,\alpha,\pi)$ is
constructed using the following iterative procedure:
\begin{figure*}[!htp]
\begin{center}
\fbox{
\begin{minipage}{\textwidth}\it
Let $A_1$ be the set of variables whose support w.r.t. $\varphi$ is at most $d/3$.\\
Let $A_2$ be the set of non-stable variables w.r.t. $\pi$.\\
Let $A_3$ be the set of stable variables w.r.t. $\pi$ which are
violated by $\alpha$.
\begin{enumerate}
\item Set $H_{0} = V \setminus (A_1 \cup A_2 \cup A_3)$.
\item While there exists a variable $a_{i} \in H_{i}$ which supports
less than $d/4$ clauses in $\F[H_{i}]$ OR appears in more than
$d/30$ clauses \emph{not} in $\F[H_{i}]$ define $H_{i+1}=H_{i}
\setminus \{a_{i}\}$.
\item Let $a_{m}$ be the last variable removed in step 2. Define
$\Core=H_{m+1}$.
\end{enumerate}\rm
\end{minipage}
}
\end{center}
\end{figure*}
\begin{prop}\label{SizeOfHBarPr}
If both $\alpha$ and $\pi$ are chosen uniformly at random then
with probability $1-n^{-\Omega(d)}$, $\#\Core = (1-e^{-\Omega(d)})n$.
\end{prop}
\begin{Proof} Let $\bar{\Core}=V\setminus \Core$. Set $\delta = e^{-\Theta(d)}$.
Partition the variables in $\bar{\Core}$ into variables that
belong to $A_1 \cup A_2 \cup A_3$, and variables that were removed
in the iterative step, $\bar{H}^{it}=H_{0}\setminus \Core$. If
$\#\bar{\Core}\geq \delta n$, then at least one of $A_1 \cup A_2
\cup A_3$, $\bar{H}^{it}$ has cardinality at least $\delta n/2$.
Consequently,
\begin{align*}
Pr[\#\bar{\Core}\geq \delta n] \leq \underbrace{Pr[\#A_1 \cup A_2
\cup A_3\geq \delta n/2]}_{(a)}+\underbrace{Pr[\#\bar{H}^{it} \geq
\delta n/2\text{ }\big|\text{ }\#A_1 \cup A_2 \cup A_3 \leq\delta
n/2]}_{(b)}.
\end{align*}
Propositions \ref{SupportSuccRate}, \ref{StableSuccRate}, and
\ref{ViolatedSuccRate} are used
to bound $(a)$. To bound $(b)$, observe that every variable that
is removed in iteration $i$ of the iterative step (step 2),
supports at least $(d/3-d/4)=d/12$ clauses in which at least
another variable belongs to $\{a_{1},a_{2},...,a_{i-1}\}\cup A_1
\cup A_2 \cup A_3$, or appears in $d/30$ clauses each containing
at least one of the latter variables. Consider iteration $\delta
n/2$. Assuming $\#A_1 \cup A_2 \cup A_3 \leq\delta n/2$, by the
end of this iteration there exists a set containing at most
$\delta n$ variables, and there are at least $d/30\cdot\delta
n/2\cdot 1/3$ clauses containing at least two variables from it
(we divide by 3 as every clause might have been counted 3 times).
Plugging $c=d/180$ in Proposition \ref{NoDenseSubformulas}, $(b)$
is bounded. Finally observe that $(a)$ occurs with exponentially small probability, and
$(b)$ occurs with probability $n^{-\Omega(d)}$.
\end{Proof}
\subsection{The Factor Graph of the non-Core Variables}\label{NonCoreStructSec}
Proposition \ref{SizeOfHBarPr} implies that for $p=c\log n/n^2$, $c$
a sufficiently large constant, $\whp$ $\Core$ contains already all
variables. Therefore the following propositions are relevant for the
setting of Theorem \ref{ConvergenceThmMed} (namely, $p=O(1/n^2)$).
In what follows we establish the typical structure of the factor graph induced on the non-core
variables.

\begin{prop}\label{prop_NonCoreFactorGraph} Let $\F$ be a 3CNF formula
randomly sampled according to $\PlantedDist$, where $p = d/n^2$,
$d \geq d_0$, $d_0$ a sufficiently large constant, let $\pi$ be the initial random ordering of
the clause-variable messages, and $\alpha$ the initial random clause-variable message vector.
Let $T$ be the factor graph induced on the non-core
variables. $T$ enjoys the following properties:
\begin{enumerate}
  \item Every connected component contains $\whp$ $O(\log n)$ variables,
  \item every connected component contains $\whp$ at most one cycle,
  \item the probability of a cycle of length at least $k$ in $T$ is at most $e^{-\Theta(dk)}$.
\end{enumerate}
\end{prop}

A proposition of similar flavor to $(a)$ was proven in
\cite{flaxman} though with respect to a different notion of core. This alone would not have sufficed
to prove our result, and we need $(b)$ and $(c)$ as well.

\begin{cor}\label{cor:NoLongCycle}  Let $f=f(n)$ be such that $f(n)\rightarrow \infty$ as $n\to \infty$. Then $\whp$ there is no cycle of length $f(n)$ in the non-core factor graph.
\end{cor}

The proof of Proposition \ref{prop_NonCoreFactorGraph} is quite long and technical.
To avoid distraction, we go ahead and prove Theorem \ref{ConvergenceThmMed} and Proposition \ref{ConvergencePropDense}, and defer the proof
of Proposition \ref{prop_NonCoreFactorGraph} to Section \ref{sec:NoTwoCycleProof}.

\section{Proof of Theorem \ref{ConvergenceThmMed} and Proposition \ref{ConvergencePropDense}}\label{ResultsProofSection}
We start by giving an outline of the proof of Theorem
\ref{ConvergenceThmMed}. Proposition \ref{ConvergencePropDense} is derived as an easy corollary of that
proof. To prove Theorem \ref{ConvergenceThmMed} we need to establish three properties:

\begin{enumerate}

\item {\em Convergence}: the WP algorithm converges to a fixed point.

\item {\em Consistency}: the partial assignment implied by this fixed point is consistent with some satisfying assignment.

\item {\em Simplicity}: after simplifying the input formula by substituting in the values of the assigned variables,
the remaining subformula is not only satisfiable (this is handled by consistency), but also simple.

\end{enumerate}

We assume that the formula $\F$ and the execution of
$\textsf{WP}$ are {\em typical} in the sense that Propositions
\ref{SizeOfHBarPr} and \ref{prop_NonCoreFactorGraph} hold. First we prove that after one iteration
\textsf{WP} sets the core variables $\Core$ correctly ($B_i$
agrees with $\varphi$ in sign) and this assignment does not change
in later iterations. The proof of this property is rather
straightforward from the definition of a core. This establishes convergence and consistency for the core variables. From
iteration 2 onwards $\WP$ is basically running on $\F$ in
which variables belonging to $\Core$ are substituted with their
planted assignment. This subformula is satisfiable. Moreover, its factor
graph contains small (logarithmic size) connected components, each
containing at most one cycle. This last fact serves a dual purpose. It shows that if WP will eventually converge, the simplicity property will necessarily hold. Moreover, it will assist us in proving convergence and consistency for the subformula. Consider a connected component
composed of a cycle and trees ``hanging" on the cycle. Proving
convergence on the trees is done using a standard inductive
argument. The more interesting part is proving convergence on the
cycle. The difficulty there is that messages on a cycle may have more than one fixed point to which they may possibly converge,
which makes it more difficult to prove that they converge at all. Our proof starts with a case analysis that identifies those cases that have
multiple fixed points. For these cases we prove that almost surely random fluctuations caused by step~3.a of the WP algorithm will lead to convergence to some fixed point. This is
similar in flavor to the fact that a random-walk on a line eventually reaches an endpoint of the line (even though one cannot tell a-priori which endpoint this will be). Hand-in-hand with establishing convergence for the trees and cycle, we shall also prove consistency.

The set $V_A$ of Theorem~\ref{ConvergenceThmMed} is composed of all
variables from $\Core$ and those variables from the non-core factor
graph that get assigned. The set $V_U$ is composed of the UNASSIGNED
variables from non-core factor graph. We now proceed with the formal
proof.

\subsection{Analysis of \textsf{WP} on the core factor graph}
We start by proving that the messages concerning the factor graph
induced by the core-variables converge to the correct value, and
remain the same until the end of the execution.

We say that a message $C \rightarrow x$, $C=(\ell_x \vee \ell_y \vee
\ell_z)$, is \emph{correct} if its value is the same as it is when
$y \rightarrow C$ and $z\rightarrow C$ are 1 (that is agree in sign with their
planted assignment). In other words, $C \rightarrow x$ is 1 iff $C=(x\vee \bar{y}\vee\bar{z})$
($x$ supports $C$).

\begin{prop} \label{prop:invariant}
If $x_i\in \Core$ and all messages $C \rightarrow x_i$, $C \in
\F[\Core]$ are correct at the beginning of an iteration (line 3 in
the \textsf{WP} algorithm), then this invariant is kept by the end
of that iteration.
\end{prop}
\begin{Proof}
By contradiction, let $C_0 \rightarrow x$ be the first wrongly
evaluated message in the iteration. W.l.o.g. assume $C_0=(\ell_x\vee
\ell_y \vee \ell_z)$. Then at least one of $y,z$ sent a wrong
message to $C_0$. $$y \rightarrow C_0 = \sum_{C \in N^+(y),C\neq
C_0} C \rightarrow y - \sum_{C'\in N^-(y),C'\neq C_0}C' \rightarrow
y.$$ Every message $C'' \rightarrow y$, $C''\in F[\Core]\cap
\{N^{++}(y)\cup N^{-}(y)\}$ is 0 (since it was correct at the
beginning of the iteration and that didn't change until evaluating
$C_0 \rightarrow x$). On the other hand, $y \in \Core$ and therefore
it supports at least $d/4$ clauses in $\F[\Core]$. Thus at least
$(d/4-1)$ messages in the left hand sum are `1' (we subtract 1 as
$y$ might support $C_0$). $y$ appears in at most $d/30$ clauses with
non-core variables (all of which may contribute a wrong '1' message
to the right hand sum). All in all, $y \rightarrow C_0 \geq
(d/4-d/30-1)
> d/5$, which is correct (recall, we assume $\varphi = \one^n$). The
same applies for $z$, contradicting our assumption.
\end{Proof}

\begin{prop}\label{PropNoProof}
If $x_i\in \Core$ and all messages $C \rightarrow x_i$, $C \in
\F[\Core]$ are correct by the end of a \textsf{WP} iteration, then
$B_i$ agrees in sign with $\varphi(x_i)$ by the end of that
iteration.
\end{prop}
\noindent Proposition \ref{PropNoProof} follows immediately from the
definition of $\Core$ and the message $B_i$. It suffices to show
then that after the first iteration all messages $C \rightarrow
x_i$, $C \in \F[\Core]$ are correct.
\begin{prop}\label{CoreSetCorrectProp} If $\F$ is a typical instance in the setting of Theorem \ref{ConvergenceThmMed},
then after one iteration of $\textsf{WP(\F)}$,  for every variable
$x_i\in \Core$, every message $C \rightarrow x_i$, $C \in \F[\Core]$
is correct.
\end{prop}
\begin{Proof}
The proof is by induction on the order of the execution in the
first iteration. Consider the first message $C \rightarrow x$,
$C=(\ell_x\vee \ell_y \vee \ell_z)$, $C \in \F[\Core]$, to be
evaluated in the first iteration. Now consider the message $y
\rightarrow C$ at the time $C \rightarrow x$ is evaluated. All
messages $C' \rightarrow y$, $C'\in\F[H]$ have their initial
random value (as $C \rightarrow x$ is the first core message to be
evaluated). Furthermore, $y \in \Core$, and therefore there are at
most $d/30$ messages of the form $C'' \rightarrow y$, $C'' \notin
\F[\Core]$. $x \in \Core$ hence it is stable w.r.t. $\pi$ and not
violated by the initial clause-variable random messages. Therefore
$$y \rightarrow C \geq \underbrace{d/7}_{\text{property (c) in defn. \ref{ViolatedDefn}}}
-\underbrace{d/30}_{\text{property $(b)$ in defn.
\ref{ViolatedDefn}}}-\underbrace{d/30}_{\text{non-core
messages}}>d/14.$$ The same applies to $z$, to show that $C
\rightarrow x$ is correct. Now consider a message $C \rightarrow x$
at position $i$, and assume all core messages up to this point were
evaluated correctly. Observe that every core message $C' \rightarrow
y$ that was evaluated already, if $C'\in \{N^{++}(y)\cup
N^{-}(y)\}\cap \F[\Core]$ then its value is '0' by the induction
hypothesis. Since $x$ is not violated by $\alpha$, property $(b)$ in
definition \ref{ViolatedDefn} ensures that to begin with
$|\#\one_\alpha(N^{++}(y))-\#\one_\alpha(N^{-}(y))|\leq d/30$. $y
\in \Core$, therefore it appears in at most $d/30$ non-core
messages, all of which could have been already wrongly evaluated,
changing the above difference by additional $d/30$. As for the core
messages of $y$ which were already evaluated, since they were
evaluated correctly, property $(a)$ in definition \ref{ViolatedDefn}
ensures that the above difference changes by at most additional
$d/30$. All in all, by the time we evaluate $C \rightarrow x$,
$$\sum_{C' \in N^{++}(y),C'\neq C} C' \rightarrow y - \sum_{C''\in
N^-(y),C''\neq C}C'' \rightarrow y \geq -3\cdot d/30.$$ As for
messages that $y$ supports, property $(c)$ in definition
\ref{ViolatedDefn} ensures that their contribution is at least $d/7$
to begin with. Every core message in $N^{s}(y)$ that was evaluated
turned to '1', every non-core message was already counted in the
above difference. Therefore $y \rightarrow C \geq d/7-3\cdot
d/30>d/25$. The same applies to $z$ showing that $C \rightarrow x$
is correct.
\end{Proof}

To prove Proposition \ref{ConvergencePropDense}, observe that when
$p=c\log n/n^2$, with $c$ a sufficiently large constant,
Proposition \ref{SizeOfHBarPr} implies $\Core = V$. Combining this
with Proposition \ref{CoreSetCorrectProp}, Proposition
\ref{ConvergencePropDense} readily follows.

\subsection{The effect of messages that already converged}

It now remains to analyze the behavior of \textsf{WP} on the
non-core factor graph, given that the messages involving the core
factor graph have converged correctly. A key observation is that once the
messages in the factor graph induced by the core variables
converged, we can think of $\textsf{WP}$ as if running on the
formula resulting from replacing every core variable with its
planted assignment and simplifying (which may result in a
1-2-3CNF). The observation is made formal by the following proposition:


\begin{prop}\label{prop:AsIfSimplified}
Consider a run of $\WP$ that has converged on the core. Starting
at some iteration after $\WP$ has converged on the core, consider
two alternative continuations of the warning propagation
algorithm. $\textsf{WP}_1$ denotes continuing with $\WP$ on the
original input formula. $\textsf{WP}_2$ denotes continuing with
$\textsf{WP}$ on the formula obtained by replacing each core
variable with its planted assignment and simplifying. Then for
every iteration $t$, the sequence of messages in the $t$'th
iteration of $\textsf{WP}_2$ is identical to the respective
subsequence in $\textsf{WP}_1$. (This subsequence includes those
messages not involving the core variables, and includes messages
of type $x \to C$ and of the type $C \to x$.)
\end{prop}

\begin{Proof}
First note that all messages $x \to C$, $x\in \Core$,
do not change (sign) from the second iteration onwards (by the
analysis in the proof of Proposition \ref{CoreSetCorrectProp}).
Furthermore, if $\ell_x$ satisfies $C$ in $\varphi$, then $x\to C$
is positive (if $x$ is a true literal in $C$, or negative
otherwise), and therefore all messages $C\to y$, $y\neq x$ are
constantly 0. Namely, they don't effect any calculation, and this
is as if we replaced $\ell_x$ with TRUE, and in the simplification
process $C$ disappeared. If $\ell_x$ is false in $C$ under
$\varphi$, then $x \to C$ is constantly negative (if $\ell_x=x$,
or constantly positive if $\ell_x=\bar{x}$), and this is exactly
like having $\ell_x$ removed from $C$ (which is the result of the
simplification process).

\end{Proof}

\subsection{Analysis of \textsf{WP} on the \emph{non}-core factor graph}
\label{sec:noncoreoutline}

Note that to prove the convergence of the algorithm we need also to prove that messages
of the sort $C \to x$ where $C$ is not in the core and $x$ is in the core converge.
However, if
we prove that all messages in the factor graph induced by the
non-core variables converge, then this (with the fact that the
core factor graph messages converge) immediately implies the
convergence of messages of this type.
Therefore, our {\em goal reduces to proving
convergence of \textsf{WP} on the factor graph induced by
$\F|_\psi$, where $\psi$ assigns the core variables their planted
assignment, and the rest are UNASSIGNED.}

We say that $\WP$ converged correctly in a connected component $\cal{C}$ of the non-core factor graph if
there exists a satisfying assignment $\psi$ of the entire formula which is consistent with the planted assignment on the core, and with the assignment of $\WP$ to $\cal{C}$.

\medskip

Consider a connected component in the non-core factor graph consisting of a
cycle with trees hanging from it. Our analysis proceeds in three steps:

\begin{enumerate}
\item
We first prove that
clause-variable and variable-clause messages of the form $\alpha
\to \beta$ where $\alpha \to \beta$ lead from the trees to the cycle,
converge weakly correctly w.r.t. the planted assignment. In the case that the
component has no cycles, this concludes the proof.
\item
Then,
using a refined case analysis, we show that the
messages along the cycle also converge $\whp$, this time not
necessarily to the planted assignment, but to some satisfying
assignment which agrees with the already converged messages.
\item
We conclude by showing that messages in the direction from the cycle to
the trees converge. Finally we show that together, all messages (including parts $(a)$ and $(b)$) in the connected
component converge correctly according to some satisfying
assignment.
\end{enumerate}

Consider the factor graph $F$ induced by the simplified formula. A
{\em cycle} in $F$ is a collection $x_1,C_2,x_3,C_4,\ldots,x_r =
x_1$ where $x_i$ and $x_{i+2}$ belong to $C_{i+1}$ for all $i$ (in
our description we consider only odd values of $i$) and $x_i \neq
x_{i+2}$, $C_{i+1} \neq C_{i+3}$ for all $i$. A factor graph $F$ is
a {\em tree} if it contains no cycles. It is {\em unicyclic} if it
contains exactly one cycle. Let $x \to C$ be a directed edge of $F$.
We say that $x \to C$ {\em belongs} to the cycle, if both $x$ and
$C$ belong to the cycle. For an edge $x \to C$ that does not belong
to the cycle, we say that $x \to C$ {\em is directed towards} the
cycle if $x$ doesn't belong to the cycle and $C$ lies on the simple path from $x$ to the cycle. We say that the edge $x
\to C$ is {\em directed away} from the cycle if $C$ doesn't belong
to the cycle and $x$ lies on the simple path from the cycle to $C$. Similarly we define what it means for an
edges $C \to x$ to belong to the cycle, to be directed towards the
cycle and to be directed away from the cycle.
\begin{prop}
Let $F$ be a unicyclic factor graph. Then every directed edge of the form $x \to C$ or $C \to x$ either belongs to the cycle, or is directed towards it or directed away from it.
\end{prop}
\begin{Proof}
Recall that the factor graph is an undirected graph, and the
direction is associated with the messages. Take an edge $x\to C$
(similarly for $C\to x$), if it lies on the cycle, then we are done.
Otherwise, since the factor graph is connected, consider the path in
the tree leading from some element of the cycle to $C$. This path is
either contained in the path to $x$ or contains it (otherwise there
is another cycle). In the first case $x\to C$ is directed towards
the cycle, and in the latter $x\to C$ is directed away from the
cycle.
\end{Proof}

Our analysis proceeds in two parts: first we shall analyze $\WP$ on the trees, then
$\WP$ on the cycle and connect the two (which is relevant for the uni-cyclic components).

\subsection{$\WP$ on the trees}
As we already mentioned before, there are two directions to
consider: messages directed towards the cycle and away from the
cycle. In this section we shall consider a rooted tree, and
partition the messages according to messages which are oriented
away from the root (they will correspond in the sequel to messages
going away from the cycle) and messages towards the root (messages
from the leaves towards the root -- later to be identified with
messages going into the cycle). The first lemma concerns messages
going towards the root.



\begin{rem}\rm\label{rem:ConvergenceOnTrees} Lemma \ref{MessageIntoCycleConvergeLem}
is a special case of the known fact (see~\cite{SurveyPropagation}
for example) that for every tree induced by a satisfiable formula,
$\WP$ converges and there exists a satisfying assignment $\psi$
such that every $B_i$ is either 0 or agrees with $\psi$. In Lemma
\ref{MessageIntoCycleConvergeLem} we assume that the formula is
satisfiable by the all~1 assignment (the planted assignment), and
consider only messages towards the root.
\end{rem}

\begin{lem}\label{MessageIntoCycleConvergeLem}
Let $C \rightarrow x$ be an edge in the non-core factor graph
belonging to a connected component of size $s$, and in particular to a rooted tree $T$. If
$C\rightarrow x$ is directed towards the root then the message $C
\rightarrow x$ converges after at most $O(s)$ iterations.
Furthermore, if $C \rightarrow x = 1$ then $x$ appears positively in
$C$.
\end{lem}
\begin{Proof}
We consider the case $C=(\ell_x \vee \ell_y)$ -- the case
$C=(\ell_x \vee \ell_y \vee \ell_z)$ where all three literals
belong to non-core variables is proved similarly. For an edge
$(C,x)$ in the factor graph, we define $\verb"level"(C,x)$ to be
the number of edges in a path between $C$ and the leaf most
distant from $C$ in the factor graph from which the edge $(C,x)$
is removed. The lemma is now proved using induction on the level
$i$. Namely, after the $i^{th}$ iteration, all messages $C
\rightarrow x$ associated with an edge $(C,x)$ at level $i$
converge, and if $C \rightarrow x = 1$ then $x$ appears positively
in $C$.

The base case is an edge $(C,x)$ at level 0. If
$\verb"level"(C,x)=0$ then $C$ is a unit clause containing only
the variable $x$. By the definition of the messages, in this case
$C\rightarrow x = 1$ and indeed it must be the case that $x$ is
positive in $C$ (as the other two variables evaluate to FALSE
under the planted). Now consider an edge $(C,x)$ at level $i$, and
consider iteration $i$. Since $i>0$, it must be that there is
another non-core variable $y$ in $C$ (or two more variables
$y,z$). Consider an edge $(C',y)$, $y\in C'$ (if no such $C'$
exists that we are done as $C \to x$ will be constantly 0 in this
case).

$\verb"level"(C',y)$ is strictly smaller than $i$ since every path
from $C$ to a leaf (when deleting the edge $(C,x)$) passes through
some edge $(C',y)$. By the induction hypothesis, all messages $C'
\rightarrow y$ already converged, and therefore also $y\to C$ and
in turn $C\to x$. It is only left to take care of the case $C\to
x=1$. In this case, there must be a clause $C'$ s.t. $C'
\rightarrow y = 1$ and $y$ appears positively in $C'$ (by the
induction hypothesis). If $C \to x=1$ it must be that $y$ appears
negatively in $C$ and therefore $x$ must appear positively
(otherwise $C$ is not satisfied by the planted assignment).
\end{Proof}


Next we consider several scenarios that correspond to messages going
from the root towards the leaves. Those scenarios correspond to
step~(c) of our analysis, referred to in
Section~\ref{sec:noncoreoutline}.

\begin{clm}\label{1MsgIffUnsatClm}
Let $F$ be a unicyclic formula, and assume that WP has converged on $F$.
Let $x \to C$ be directed away from the cycle, and define $F_C$ to be the
subformula inducing the tree rooted at $C$ ($x$ itself is removed from $C$).
Then $C \to x = 0$ in the fixed point if and only if $F_C$ is satisfiable.
\end{clm}
\begin{Proof}
The structure of the proof is similar to that of Lemma
\ref{MessageIntoCycleConvergeLem}.
For convenience we extend the definition of level above as to include edges on the cycle.
We say that an edge $(C,x)$
in the factor graph has $\verb"level"(C,x)$ equal $\infty$ if
$(C,x)$ lies on a cycle and  $\verb"level"(C,x) = t < \infty$ if $t$ is the maximal length of a
path between $C$ and a leaf in the factor graph from which the edge
$(C,x)$ is removed. The lemma is now proved using induction on $t$.

The base case is an edge $(C,x)$ at level 0. If
$\verb"level"(C,x)=0$ then $C$ is a unit clause containing only the
variable $x$, and then $F_C$ is the empty formula. Indeed $C \to x =1$ by definition and
$F_C$ is unsatisfiable (by definition again, the empty formula is not satisfiable).

Now consider an edge $(C,x)$ at level $t>0$ ($t$ is still finite because we are not
interested in cycle edges). Assume w.l.o.g. that $C=(x\vee
y \vee z)$ (maybe only $y$). Let$C_1,C_2,\ldots,C_r$ be the clauses in $F_C$ that contain the variable $y$,
and let $D_1,D_2,\ldots,D_s$ be the clauses that contain the variable $z$. Similarly to $F_C$ we can define $F_{C_i}$
and $F_{D_j}$. Observe that every edge
$(C_i,y)$ and $(D_j,z)$ in $F_C$ has $\verb"level"$ value which is strictly
smaller than $t$ -- since every path from $C$ to a leaf (when
deleting the edge $(C,x)$) passes through some such edge.

First consider the case where $C \to x = 1$. For that to happen, it must be that
there are two clauses $C_i,D_j$ in $F_C$ so that $C_i \to y = 1$ and $D_j \to z = 1$ and
$C_i$ contains $\bar{y}$, $D_j$ contains $\bar{z}$. By the induction hypothesis, $F_{C_i}$ and
$F_{D_j}$ are not satisfiable. The only way that $F_C$ can now be possibly satisfied is by assigning $z$ and $y$ to FALSE,
but then $F_C$ is not satisfied as the clause $C$ (without $x$) evaluates to FALSE.

Next we prove that if $F_C$ is not satisfiable then it must be that $C\to
x=1$. If $F_C$ is not satisfiable then it must be that
$C=(\ell_x\vee \bar{y} \vee \bar{z})$ (otherwise $\varphi$ satisfies $F_C$). Further, observe that
there must exist at least one clause $C_i$ containing $y$
positively, and at least one $D_j$ containing $z$ positively s.t.
$F_{C_{i}}$ and $F_{D_{j}}$ are unsatisfiable. Otherwise, we can
define $\varphi'$ to be $\varphi$ except that $y$ or $z$ are
assigned FALSE (depending which of $C_i$ or $D_j$ doesn't exist). It
is easy to see that $\varphi'$ satisfies $F_C$, contradicting our
assumption. By the induction hypothesis $C_i\to y=1$ and $D_j\to
z=1$. Further, by Lemma \ref{MessageIntoCycleConvergeLem}, there cannot be a message
$C' \to y = 1$ where $y$ is negative, or $D' \to z =1$ where $z$ is negative (because all of these messages are directed towards the cycle). This in turn implies that $C\to x=1$.
\end{Proof}\\
Recall our definition for $\WP$ converges correctly on a subformula $F'$ of $F$ if
there exists a satisfying assignment $\psi$ to $F$ which is consistent with the planted assignment $\varphi$ on the core,
and with the assignment of $\WP$ to $F'$.
\begin{clm}\label{spinalDecomp}
Assume that $F$ is a unicyclic formula. Assume further that WP has converged on $F$.
Let $C \to y$ be directed away from the cycle.
Consider a subformula $F_C$ which induces a tree rooted at a clause $C$. This formula contains the clause $C$ and
all other clauses whose path to the cycle goes via $y$.
If in the fixed point for $F$ it holds that
\begin{itemize}
  \item $C\to y=1$,
  \item $y$ appears negatively in $C$,
  \item $y\to C \geq 1$,
\end{itemize}
then $\WP$ converges correctly on $F_C$.
\end{clm}
\begin{Proof}
We split the variables of $F_C$ into two sets. The first is the set of spinal variables (to be defined shortly), and the remaining variables.
The \emph{spine} (rooted at $y$) of a tree-like formula is defined using the following procedure:
Let us track the origin of the message $y\to C$ in the tree
(which is directed towards the root). For $y\to C \geq 1$ to occur, there must be a clause $D_1$ in the tree that contains $y$ positively
and a message $D_1  \to y = 1$ in the direction of the root
(as messages in the direction of the root are only effected by other messages in that direction). Let us backtrack one more step.
 $D_1=(y \vee \bar{k} \vee \bar{w})$ for some variables $w,k$; $k$ and $w$ must appear negatively in $D_1$ by Lemma \ref{MessageIntoCycleConvergeLem},
 and the fact that $D_1 \to y = 1$, that is both $k$ and $w$ were issued warnings having them not satisfy $D_1$. Let us consider the clauses $D_2$ and $D_3$ that issues warnings to $k$ and $w$ respectively.
$D_2=(k\vee ...),D_3=(w\vee ...)$, $D_2 \to k =1$ and $D_3 \to w = 1$, and both messages are directed towards the root.
Obviously, one can inductively continue this backtracking procedure which terminates at the leaves of the tree
(since there are no cycles the procedure is well defined and always terminates).
Let us call the clauses and variables that emerge in this backtrack the {\em spine} of the tree.
The figure below illustrates this procedure.

\begin{figure*}[h]

\begin{center}
\includegraphics[width=1.7in]{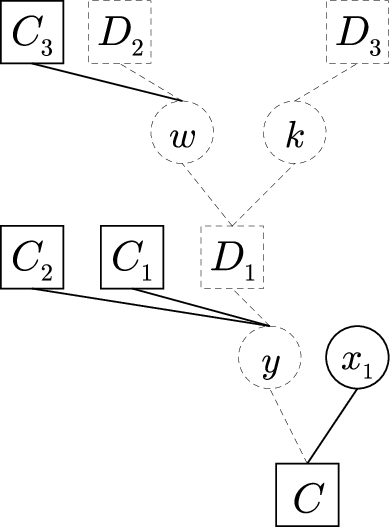}
\end{center}

\caption{The spine of a tree}
\label{fig:a}
\end{figure*}

Let us start with the $B_i$-messages for the non-spinal variables (the $B_i$ message was defined in Section \ref{WPSubs}). The spine splits $F_C$ into sub-trees hanging from the spinal variables (the root of each such tree is a clause, e.g. $C_1,C_2,C_3$ in our example). Every non-spinal variable belongs to one of those trees. It suffices to prove that for every subtree rooted at $C_i$, hanging from a spinal variable $s$, it holds that $s \to C_i \geq 0$ in the fixed point (that is, $s$ transmits its assignment according to the planted assignment). This ensures, by Remarks \ref{rem:ConvergenceOnTrees}, convergence for that tree to a satisfying assignment (consistent with the planted). Let us prove that $s \to C_i \geq 0$. Let $F_s$ be the subformula corresponding to the tree hanging from $s$ rooted at some clause $C_i$. We prove that $s \to C_i \geq 0$ by induction on
the distance in $F_C$ between $y$ and $s$. The base case is distance 0, which is $y$ itself. The messages that we
need to verify are of the form $y \to C_i$, $C_i \neq D_1$, which are
pointing away from the cycle. For every message $y \to C_i$, the wrong message
$C \to y$ is evened by the correct warning $D_1 \to y$. Since $y\to C_i$
depends only on one message which is directed away from the cycle (or else
 there is a second cycle), and all the other messages that contribute to the calculation of $y\to C_i$ are correct (Lemma \ref{MessageIntoCycleConvergeLem}), we conclude that $y \to C_i \geq 0$. The induction step follows very
similarly.

Let us now consider the spinal variables. For every such variable there is always at most one message in the direction away
from the cycle (otherwise there is more than one cycle); this
message may be wrong. Using a very similar inductive argument as in the previous paragraph, one
can show that there is always at least one correct warning (in the
direction of the cycle) for each spinal variable. Therefore $B_s \geq 0$ for every spinal
variable $s$.
\end{Proof}

\subsection{$\WP$ on cycles} \label{subsec:cycle}

We will denote a cycle by $x_1, C_2 , x_3, C_4
...x_{2r-1},C_{2r},x_1$ where by this we mean that $x_i$ appears
in the clauses before/after it and that $C_i$ contains the two
variables before/after it. We consider two different types of
cycles.
\begin{itemize}
\item
\emph{Biased} cycles: cycles that have at least one warning
message $C \to x_i = 1$ coming into the cycle, where $C \to x_i$ directs into the cycle
and the value of $C \to x_i$ is the value after the edge has converged.
\item
\emph{Free} cycles: cycles that do not have such messages coming
in, or all messages coming in are 0 messages.
\end{itemize}

\subsubsection{Convergence of \textsf{WP} when the cycle is
biased}

\noindent First we observe that we may assume w.l.o.g. that edges that enter the cycle enter it at a variable rather than at a clause (hence that every
clause on the cycle contains exactly two non-core variables). This is
because of a simple argument similar to Proposition
\ref{prop:AsIfSimplified}: consider an edge going into the cycle,
$z\to C$, and w.l.o.g. assume that $z$ appears positively in $C$.
After all the edges going into the cycle have converged, if $z\to
C\geq 0$ it follows that $C\to x=0$ for cycle edges $(C,x)$, and
thus execution on the cycle is the same as if $C$ was removed from
the formula, only now we are left with a tree, for which convergence
to a correct assignment is guaranteed (Remark
\ref{rem:ConvergenceOnTrees}). If $z\to C<0$, then the execution is
exactly as if $z$ was removed from $C$ (and $C$ is in 2-CNF form).

\begin{prop}\label{prop:ConnvergenceInBiasedCaseClm}
Let $\cal{C}$ be a connected component of the factor graph of size
$s$ containing one cycle s.t. there exists an edge directed into the
cycle $C \to x_i$ where $x_i$ belongs to the cycle and such that the
message converges to $C \to x_i = 1$. Then $\WP$ converges on
$\cal{C}$ after at most $O(s)$ rounds. Moreover for the fixed point,
if the message $C' \to x = 1$ then $x$ appears positively in $C'$.
\end{prop}

\begin{Proof}
A message of the cycle $C_j \to x_{j+1}$
depends only on cycle messages of the type $C_{j'} \to x_{j'+1},
x_{j'+1} \to C_{j'+2}$ and on messages coming into the cycle. In
other words during the execution of $\WP$ the values of all messages
$C_{j'} \to x_{j'-1},x_{j'-1} \to C_{j'-2}$ do not effect the value
of the message $C_j \to x_{j+1}$. Recall that we are in the case
where there exists a message $C \to x_i = 1$ going into the cycle
(after the convergence of these messages). Also $x_i$ must
appear positively in $C$ (Lemma \ref{MessageIntoCycleConvergeLem}). We consider the following cases:

\begin{itemize}
\item
There exists a variable $x_j$ that appears positively in both $C_{j-1}$ and $C_{j+1}$ (the case $j=i$ is allowed here).
We note that in this case the message $x_j \to C_{j+1}$ must be non-negative which implies that the message $C_{j+1} \to x_{j+2}$ converges to the value $0$. This in turn implies that the value of all messages $x_r \to C_{r+1}$
and $C_{r+1} \to x_{r+2}$ for $r \neq j$ will remain the same if the clause
$C_{j+1}$ is removed from the formula.
However, this case reduces to the case of tree formula.
\item
$x_i$ appears negatively in $C_{i+1}$ and positively in $C_{i-1}$. We note that in this case the value of the message $x_i \to C_{i+1}$ is always at least 1, which implies that the message $C_{i+1} \to x_{i+2}$ always take the value $1$. Thus in this case we may remove the clause $C_{i+1}$ from the formula and replace it by the unit clause $\ell_y$ where $C_{i+1} = \ell_y \vee \bar{x_i}$. Again, this reduces to the case of a tree formula.
\item
The remaining case is the case where $x_i$ appears negatively in both $C_{i-1}$ and $C_{i+1}$ and there is no $j$ such that $x_j$ appears positively in both $C_{j-1}$ and $C_{j+1}$. We claim that this leads to contradiction. An easy corollary of Lemma \ref{MessageIntoCycleConvergeLem} is that all the messages that go into the cycle have converged according to the planted assignment. Therefore w.l.o.g one can ignore messages of the form $y \to C_i$, $y$ is a non-cycle variable, and $C_i$ is a cycle clause (we can assume that every such message says that the literal of $y$ is going to be false in $C_i$, otherwise the cycle structure is again broken and we reduce to the tree case). Therefore the cycle reduces to a 2SAT formula which has to be satisfiable by the planted assignment, in which in particular $x_i = 1$.
Write $C_{i+1} = \bar{x_i} \vee \ell_{i+2}$. Then for the satisfying assignment we must have $\ell_{i+2} = 1$, similarly
$\ell_{i+4} = 1$, etc, until we reach $\bar{x}_i = 1$, a contradiction.
\end{itemize}
To summarize, by Lemma \ref{MessageIntoCycleConvergeLem} the messages going into the rooted tree at $x_i$ converge after $O(s)$ steps, and at least one warning is issued. By the above discussion, for every clause $D$ in the connected component it holds that $x_i \to D\geq 0$ (as $x_i$ appears in at most one message which may be wrong -- a cycle message). Since there is always a satisfying assignment consistent with $x_i$ assigned TRUE, then after reducing the cycle to a tree we are left with a satisfiable tree. Remark \ref{rem:ConvergenceOnTrees} guarantees convergence in additional $O(s)$ iterations.
\end{Proof}

\subsubsection {Convergence of \textsf{WP} when the cycle is free}
The main result of this subsection is summarized in the following claim:
\begin{clm}\label{ConnvergenceInNonBiasedCaseMainClm}
Let $\cal{C}$ be a connected component of the factor graph of size
$s$ containing one cycle of size $r$ s.t. the fixed point contains
no messages $C \to x = 1$ going into the cycle (the cycle is free).
Then $\whp$ \textsf{WP} converges on $\cal{C}$ after at most $O(r^2 \cdot \log
n+s)$ rounds. Moreover for the fixed point, if we simplify the
formula which induces $\cal{C}$ according to the resulting $B_i$'s,
then the resulting subformula is satisfiable.
\end{clm}

\begin{rem}\label{DiffRem}\rm Observe that the free case is the only one where
convergence according to the planted assignment is not guaranteed.
Furthermore, the free cycle case is the one that may not converge
``quickly" (or not at all), though this is extremely unlikely. The
proof of Proposition \ref{ConnvergenceInNonBiasedCaseMainClm} is the
only place in the analysis where we use the fact that in line
$3.a$ of \textsf{WP} we use fresh randomness in every
iteration.
\end{rem}

We consider two cases: the
easy one is the case in which the cycle contains a pure variable
w.r.t the cycle (though this variable may not be pure w.r.t to the
entire formula).

\begin{clm}\label{ConnvergenceOfNonBiasedCyclePure}
If the cycle contains a variable $x_i$ appearing in the same
polarity in both $C_{i+1},C_{i-1}$, then the messages $C \to x$
along cycle edges converge. Moreover for the fixed point, if $C
\to x = 1$ then $x$ satisfies $C$ according to $\varphi$.
\end{clm}
The proof is essentially the first case in the proof of
Proposition \ref{prop:ConnvergenceInBiasedCaseClm}. We omit the details.

\noindent We now move to the harder case, in which the cycle
contains no pure variables (which is the case referred to in Remark
\ref{DiffRem}).
\begin{prop}\label{FreeCycleNoPureConvergenceClm}
Consider a free cycle of size $r$ with no pure literal, and one of
the two directed cycles of messages. Then the messages along the
cycle converge $\whp$ to either all $0$ or all $1$ in $O(r^2\log n)$
rounds.
\end{prop}
Convergence in $O(r^2\log n)$ rounds suffices due to Corollary
\ref{cor:NoLongCycle} (which asserts that $\whp$ the length of every cycle is constant).
The proof of Proposition \ref{FreeCycleNoPureConvergenceClm} is
given at the end of this section. We proceed by analyzing
\textsf{WP} assuming that Proposition
\ref{FreeCycleNoPureConvergenceClm} holds, which is the case $\whp$.

\begin{prop}\label{prop:free_cycle_comp_convergence} Suppose that
the cycle messages have converged (in the setting of Proposition
\ref{FreeCycleNoPureConvergenceClm}), then the formula resulting
from substituting every $x_i$ with the value assigned to it by $B_i$
(according to the fixed point of $\WP$), and simplifying, is
satisfiable.
\end{prop}
\begin{Proof}
Let $F$ be the subformula that induces the connected component
$\cal{C}$, and decompose it according to the trees that hang on the
cycle's variables and the trees that hang on the cycle's clauses.
Observe that the formulas that induce these trees are variable and
clause disjoint (since there is only one cycle in $\cal{C}$).

Let us start with the cycle clauses. The key observation is that setting the
cycle variables according to one arbitrary orientation (say, set
$x_i$ to satisfy $C_{i+1}$) satisfies the cycle and doesn't conflict
with any satisfying assignment of the hanging trees: if the tree
hangs on a variable $x_i$, then since the cycle is free, the tree is
satisfiable regardless of the assignment of $x_i$ (Proposition
\ref{1MsgIffUnsatClm}). In the case that the tree hangs on a
cycle-clause $C$, then the cycle variables and the tree variables
are disjoint, and $C$ is satisfied already by a cycle-variable
regardless of the assignment of the tree-variables. Now how does this coincide
with the result of $\WP$. Recall that we are in the case where the cycle is free. Therefore only messages
$C \to x_i$ where both $C$ and $x_i$ belong to the cycle effect $B_i$. If in the fixed
point one cycle orientation is 0 and one orientation is 1, then the
$B_i$ messages of the cycle variables implement exactly this policy.
If both cycle orientations converged to 1 or to 0, then the
corresponding $B_i$ messages of all cycle variables are UNASSIGNED
(since the cycle is free), but then the same policy can be used to
satisfy the clauses of the cycle in a manner consistent with the
rest of the formula.

It remains to show that $\WP$ converges on every tree in a manner that is consistent with some
satisfying assignment of the tree. We consider several cases.

Consider a tree hanging on a cycle variable $x_i$. Let $C$ be some
non-cycle clause that contains $x_i$, and $F_C$ the subformula that
induces the tree rooted at $C$. Observe that once the cycle has
converged, then the message $x_i\to C$ does not change anymore. If
$x_i \to C$ agrees with $\varphi$ there are two possibilities. Either $x_i$ satisfies $C$ under $\varphi$,
in which case $C$ always sends 0 to $F_C$, and then $\WP$ executes on $F_C$ as if $C$ is removed. Remark \ref{rem:ConvergenceOnTrees} guarantees correct convergence (as $F_C\setminus C$ is satisfiable), and as for $C$, $B_i \geq 0$ and we can set $x_i$ to TRUE so that it satisfies $C$ and is consistent with the assignment of the cycle ($B_i \geq 0$ since $x_i \geq 0$ and $C \to x_i = 0$ as we are in the free cycle case). If $x_i$ appears negatively in $C$, then $\WP$ executes as if $x_i$ was deleted from $C$. Still $F_C$ is satisfiable and correct convergence is guaranteed.

Now consider the case where $x_i \to C$ disagrees
with $\varphi$. Recall that we assume $\varphi(x_i)=TRUE$, and
therefore $x_i \to C$ is negative in the fixed point. If $x_i$ appears
negatively in $C$ then $C\rightarrow y=0$ for every $y\in C$ (since
$x_i$ signals $C$ that it satisfies it), and therefore $C$ doesn't
effect any calculation from this point onwards, and the correct
convergence of $F_C$ is again guaranteed by Remark
\ref{rem:ConvergenceOnTrees} on the convergence for satisfiable
trees. The more intricate case is if $C$ contains $x_i$ positively.
Since we are in the free case, it must hold that $C \to x=0$.
Therefore using Proposition \ref{1MsgIffUnsatClm} one obtains
that $F_C$ is satisfiable (regardless of the assignment of $x_i$), and $\WP$ will converge
as required (again Remark \ref{rem:ConvergenceOnTrees}).

Now consider a tree hanging on a cycle clause. Namely,
$C_{i+1}=(x_i\vee x_{i+2} \vee y)$, where $x_i,x_{i+2}$ are cycle
variables, and $(C_{i+1},y)$ is a tree edge. If one of the cycle
orientations converged to 0, then $C_{i+1}\to y$ converges to 0, and
then Remark \ref{rem:ConvergenceOnTrees} guarantees correct
convergence. The same applies to the case where $C_{i+1} \to y$ converges to 1 and
$y$ is positive in $C_{i+1}$ (since then we can remove $y$ from the tree, use Remark \ref{rem:ConvergenceOnTrees}
for the remaining part, then add back $y$, and set it to TRUE without causing any conflict with the tree assignment, but satisfying $C_{i+1}$ according to the planted assignment).

The delicate case remains when $C_{i+1} \to y$ converges to 1 but
$y$'s polarity in $C_{i+1}$ disagrees with $\varphi$, that is, $y$
is negative in $C_{i+1}$. The key observation is that the message
$y\to C_{i+1}$ (which is directed towards the cycle) must have
converged to a positive value (otherwise, $C_{i+1}\to x_i$ and
$C_{i+1}\to x_{i+2}$ would have converged to 0). However this complies with the scenario of Proposition
\ref{spinalDecomp}, and again correct convergence is guaranteed.
\end{Proof}

In Theorem \ref{ConvergenceThmMed} the unassigned
variables are required to induce a ``simple" formula,
which is satisfiable in linear time. Observe that the factor
graph induced by the UNASSIGNED variables consists of connected
components whose structure is a cycle with trees hanging on it, or
just a tree. A formula whose factor graph is a tree can be
satisfied in linear time by starting with the leaves (which are
determined uniquely in case that the leaf is a clause -- namely, a
unit clause, or if the leaf is a variable then it appears only in
one clause, and can be immediately assigned) and proceeding
recursively. Regarding the cycle, consider an arbitrary variable
$x$ on the cycle. By assigning $x$ and simplifying accordingly, we
remain with a tree. Since there are only two ways to assign $x$,
the whole procedures is linear in the size of the connected
component. This completes the proof of Theorem \ref{ConvergenceThmMed}.

\subsubsection{Proof of Proposition \ref{FreeCycleNoPureConvergenceClm}}
Since the cycle has no pure literal it must be of the following
form: $C_1 = (\ell_{x_1} \vee \overline{\ell}_{x_2}), C_2 =
(\ell_{x_2} \vee \overline{\ell}_{x_3}), \ldots, C_L = (\ell_{x_L}
\vee \overline{\ell}_{x_1})$ (recall the definition of cycles at the beginning of subsection \ref{subsec:cycle}).

Consider one of the directed cycles, say: $x_1 \to C_1 \to x_2 \to
\cdots$ and note that when the message $x_i \to C_i$ is updated it
obtains the current value of $C_{i-1} \to x_i$ and when the
message $C_i \to x_{i+1}$ is updated, it obtains the current value
of $x_i \to C_i$.

It thus suffices to show that the process above converges to all
$0$ or all $1$ in time polynomial in the cycle length. This we
prove in the lemma below.

\begin{clm} \label{clm:perm_process}
Consider the process $(\Gamma^i : i \geq 0)$ taking values in $\{0,1\}^{L}$.
The process is a Markov process started at $\Gamma^0 \in \{0,1\}^{L}$.

Given the state of the process $\Gamma^i$ at step $i$, the distribution of $\Gamma^{i+1}$
at round $i+1$ is defined by picking a permutation $\sigma \in
S_L$ uniformly at random and independently. Then $\Gamma^{i+1}$ is defined
via the following process:
let $\Delta_0 = \Gamma^i$ and for $1 \leq j \leq L$:
\begin{itemize}
\item
Let $\Delta_j$ be obtained from $\Delta_{j-1}$
by setting $\Delta_j(\sigma(j-1)) = \Delta_{j-1}((\sigma(j-1)+1) \mod L)$ and
$\Delta_j(r) = \Delta_{j-1}(r)$ for all $r \neq \sigma(j-1)$.
\end{itemize}
Then set $\Gamma^{i+1} = \Delta_L$.

Let $T$ be the stopping time where the process hits the state all
$0$ or all $1$. Then for all $\Gamma^0$:
\begin{equation} \label{eq:tail_perm}
Pr[T \geq 4 a L^2] \leq L 2^{-a}.
\end{equation}
for all $a \geq 1$ integer.
\end{clm}

The proof of Proposition~\ref{clm:perm_process} is based on the
fact that the process defined in this lemma is a martingale. This
is established in the following two lemmas.

\begin{lem} \label{lem:perm_process_mart}
Consider the following variant $\tilde{\Gamma^i}$ of the process
$\Gamma^i$ defined in Proposition~\ref{clm:perm_process}. In the
variant, different intervals of $0/1$ are assigned different colors
and the copying procedure is as above. Fix one color and
let $X_i$ denote the number of elements of the cycle of that color in
 $\tilde{\Gamma}^i$. Then $X_i$ is a
martingale with respect to the filtration defined by
$\tilde{\Gamma}^i$.
\end{lem}

\begin{Proof}
From the definition
of the copying process it is easy to see that
\[
\mathbb{E}[X_{i+1} | \tilde{\Gamma}^i,\tilde{\Gamma}^{i-1},\ldots] =
\mathbb{E}[X_{i+1} | X_i].
\]
We will show below that $\mathbb{E}[X_{i+1} | X_i] = 0$ and thus that
$X_{i}$ is a martingale with respect to the filtration
$\tilde{\Gamma}^i$.

Assume that $X_i = k$. Then w.l.o.g. we may assume that
 the configuration $\tilde{\Gamma}^i$
consists of an interval of $1$'s of length $k$ and an interval of
$0$'s of length $L-k$. Indeed if the configuration consists of a number of intervals of
length $Y_{i,1},\ldots,Y_{i,r}$ where $X_i = \sum_{t} Y_{i,t}$ then we may think of the different sub-intervals
as having different colors. Then proving that each of the $Y_{i,t}$
is a martingale implies that $X_i$ is an interval as needed.

We calculate separately the expected shifts in the locations of
left end-points of the $0$ and $1$ interval respectively. We
denote the two shift random variables by $L_0$ and $L_1$. Clearly
$L_0 = I_{0,1} + I_{0,2} + \ldots + I_{0,k-1}$ where $I_{0,j}$ is
the indicator of the event that the $0$ left-end point shifted by at
least $j$ and similarly for $L_1$. Note that
\[
\mathbb{E}[I_{0,j}] = \frac{1}{j!} - \frac{1}{(L-k+j)!}
\]
and that
\[
\mathbb{E}[I_{1,j}] = \frac{1}{j!} - \frac{1}{(k+j)!}.
\]
The last equation follows from the fact that in order for the $1$
interval to extend by at least $j$, the $j$ copying has to take
place in the correct order and it is forbidden that they all took
place in the right order and the interval has become a $0$
interval. The previous equation is derived similarly. Thus
\[
\mathbb{E}[(X_{i+1} - X_i) | X_i] = \mathbb{E}[L_1] - \mathbb{E}[L_0] = \sum_{j=1}^{L-k}
\left(\frac{1}{j!} - \frac{1}{(k+j)!}\right) - \sum_{j=1}^{k}
\left(\frac{1}{j!} - \frac{1}{(L-k+j)!}\right) = 0
\]
This concludes the proof that $X_i$ is a martingale. The proof
follows.
\end{Proof}\\\\
The proof of Proposition~\ref{clm:perm_process} follows by a union
bound from the following lemma where the union is taken over all
intervals.

\begin{lem} \label{lem:perm_process_tail}
Consider the process $\tilde{\Gamma^i}$ defined in
Lemma~\ref{lem:perm_process_mart}. Fix one interval and let $X_i$
denote its length. Let $T$ be the stopping time where $ X_i$
equals either $0$ or $L$. Then
\[
Pr[T \geq 4 a L^2] \leq 2^{-a}.
\]
\end{lem}

\begin{Proof}
In order to bound the hitting probability of $0$ and $L$, we need
some bounds on the variance of the martingale differences. In
particular, we claim that unless $k=0$ or $k=L$ it holds that
\[
\mathbb{E}[(X_{t+1}-X_t)^2 | X_t = k] \geq 1/2.
\]
If $k=1$ or $k=L-1$ this follows since with probability at least
$1/2$ the end of the interval will be hit. Otherwise, it is easy
to see that the probability that $X_{t+1}-X_t$ is at least $1$ is
at least $1/4$ and similarly the probability that it is at most
$-1$ is at least $1/4$. This can be verified by considering the
event that one end-point moves by at least $2$ and the other one
by at most $1$.

Let $T$ be the stopping time when $X_T$ hits $0$ or $n$. Then by a
Wald kind of calculation we obtain:
\begin{eqnarray*}
L^2 &\geq& \mathbb{E}[(X_T-X_0)^2] = \mathbb{E}[ (\sum_{t=1}^{\infty} 1(T \geq t)
(X_t-X_{t-1}))^2] \\ &=& \mathbb{E}[\sum_{t,s=1}^{\infty}
(X_t-X_{t-1})(X_s-X_{s-1}) 1(T \geq \max{t,s})] \\ &=&
\mathbb{E}[\sum_{t=1}^{\infty} (X_t-X_{t-1})^2 1(T \geq t)] \geq
\frac{1}{2} \sum_{t=1}^{\infty} P[T \geq t] = \mathbb{E}[T]/2,
\end{eqnarray*}
where the first equality in the last line follows from the fact that
if $s < t$ say then:
\[
\mathbb{E}[(X_t-X_{t-1})(X_s-X_{s-1}) 1(T \geq \max{t,s})]
\]
\[
= \mathbb{E}[(X_t-X_{t-1})(X_s-X_{s-1}) (1 - 1(T < t))]
\]
\[
= \mathbb{E}[\mathbb{E}[(X_t - X_{t-1}) (X_s - X_{s-1}) (1 - 1(T < t)) |
X_1,\ldots,X_{t-1}]]
\]
\[
=
\mathbb{E}[(X_s - X_{s-1}) (1 - 1(T < t)) E [X_t - X_{t-1} |
X_1,\ldots,X_{t-1}]] = 0.
\]

We thus obtain that $\mathbb{E}[T] \leq 2 L^2$. This implies in turn that
$Pr[T \geq 4 L^2] \leq 1/2$ and that $P[T \geq 4 a L^2] \leq
2^{-a}$ for $a \geq 1$ since $X_t$ is a Markov chain. The proposition
follows.

\end{Proof}

\section{Proof of Proposition \ref{prop_NonCoreFactorGraph}}\label{sec:NoTwoCycleProof}
We shall start with the proof of Item $(b)$. Then we shall remark how to adjust this proof to prove $(c)$.
Item $(a)$ was proven, though with respect to a slightly different definition of a core, in \cite{flaxman}. The proof
of Item $(a)$ is identical to that in \cite{flaxman}, with the adjustments made along the proof of Item $(b)$.
Details of the proof of $(a)$ are omitted.

In order to prove Proposition \ref{prop_NonCoreFactorGraph} $(b)$ it suffices to
prove that $\whp$ there are no two cycles with a simple path
(maybe of length 0) connecting the two. To this end, we consider
all possible constellations of such prohibited subgraphs and prove
the proposition using a union bound over all of them.

Every simple $2k$-cycle in the factor graph consists of $k$
variables, w.l.o.g. say $x_{1},...,x_{k}$ (all different), and $k$
clauses $C_{1},...,C_{k}$, s.t. $x_{i},x_{i+1}\in C_{i}$. The
cycle itself consists of $2k$ edges.

As for paths, we have 3 different types of paths: paths connecting
a clause in one cycle with a variable in the other (type 1), paths
connecting two clauses (type 2), and paths connecting two
variables (type 3). Clause-variable paths are always of odd
length, and clause-clause, variable-variable paths are always of
even length. A $k$-path $P$ consists of $k$ edges. If it is a
clause-variable path, it consists of $(k-1)/2$ clauses and the
same number of variables. If it is a variable-variable path, it
consists of $k/2-1$ variables and $k/2$ clauses and symmetrically
for the clause-clause path (we don't take into account the
clauses/variables that participate in the cycle, only the ones
belonging exclusively to the path).

Our prohibited graphs consist of two cycles $C_1,C_2$ and a simple
path $P$ connecting them. We call a graph containing exactly two
simple cycles and a simple path connecting them a \emph{bi-cycle}.
The path $P$ can be of either one of the three types described above. Similarly to the
bi-cycle case, one can have a cycle $C$ and a chord $P$ in it. We
call such a cycle a \emph{chord-cycle}. For parameters $i,j,k\in
[1,n]$, and $t\in \{1,2,3\}$, we denote by $B_{2i,2j,k,t}$ a
bi-cycle consisting of a $2i$-cycle connected by a $k$-path of type
$t$ to a $2j$-cycle. Similarly, we denote by $B_{2i,k,t}$ a
chord-cycle consisting of a $2i$-cycle with a $k$-path of type $t$ as
a chord.

Our goal is then to prove that $\whp$ the graph induced by the
non-core variables contains no bi-cycles and no chord-cycles.

For a fixed factor graph $H$ we let $F_{H}\subseteq \F$ be a fixed
minimal set of clauses inducing $H$, and $V(H)$ be the set of
variables in $H$. In order for a fixed graph $H$ to belong to the
factor graph induced by the non-core variables it must be that
there exists some $F_H$ s.t. $F_{H} \subseteq \F$ and that
$V(H)\subseteq \bar{\Core}$ (put differently, $V(H)\cap
\Core=\emptyset$).

Let $B=B_{2i,2j,k,t}$ (or $B=B_{2i,k,t}$ if $B$ is a chord-cycle) be
a fixed bi-cycle and $F_B$ a fixed minimal-set of clauses inducing
$B$. We start by bounding $Pr[F_{B} \subseteq \F \text{ and }
V(B)\cap \Core = \emptyset]$ and then use the union bound over all
possible bi-cycles (chord-cycles) and inducing minimal sets of
clauses. As the two events -- $\{F_B \subseteq \F\}$ and $\{V(B)\cap
\Core = \emptyset\}$ --  are not independent, the calculations are
more involved. Loosely speaking, to circumvent the dependency issue,
one needs to defuse the effect that the event $\{F_B \subseteq \F\}$
might have on $\Core$. To this end we introduce a set $\Core^*$,
defined very similarly to $\Core$ only ``cushioned" in some sense to
overcome the dependency issues (the ``cushioning" depends on $F_B$).
This is done using similar techniques to
\cite{AlonKahale97,flaxman}.

We start by defining the new set of core variables $\Core^*$ (again
w.r.t. an ordering $\pi$ of the clause-variable messages and an
initial values vector $\alpha$). The changes compared to $\Core$ are
highlighted in bold.
\begin{figure*}[!htp]
\begin{center}
\fbox{
\begin{minipage}{\textwidth}\it
Let $B_1$ be the set of variables whose support w.r.t. $\varphi$ is at most $d/3$.\\
Let $B_2$ be the set of non-stable variables w.r.t. $\pi$ where we
redefine the gap in Definition \ref{StableDefn} to be
\hspace*{15pt}$\textbf{(d/30-6)}$ in $(a)$ and $(b)$.\\
Let $B_3$ be the set of stable variables w.r.t. $\pi$ which are
violated by $\alpha$ where we redefine the gap in
\hspace*{15pt}Definition \ref{ViolatedDefn} to be
$\textbf{(d/30-6)}$ in $(a)$ and $(b)$.\\
Let $\textbf{J}\subseteq V(B)$ be the set of variables appearing
in no more than 6 different clauses in $F_B$
\begin{enumerate}
\item Set $H'_{0} = V \setminus (B_1 \cup B_2 \cup B_3 \cup (V(F_B)\setminus J))$.
\item While there exists a variable $a_{i} \in H'_{i}$ which supports
less than $d/4$ clauses in $\F[H'_{i}]$ OR appears in more than
$\textbf{(d/30-6)}$ clauses \emph{not} in $\F[H'_{i}]$, define
$H'_{i+1}=H'_{i} \setminus \{a_{i}\}$.
\item Let $a_{m}$ be the last variable removed at step 2. Define
$\Core^*=H'_{m+1}.=H'_{m+1}$.
\end{enumerate}\rm
\end{minipage}
}
\end{center}
\end{figure*}\\
Propositions \ref{StableSuccRate} and \ref{ViolatedSuccRate} could
be easily adjusted to accommodate the 6-gap in the new definition
in $B_2$ and $B_3$. Therefore Proposition \ref{SizeOfHBarPr} can
be safely restated in the context of $\Core^*$:
\begin{prop}\label{SizeOfHBarTagPr}
If both $\alpha$ and $\pi$ are chosen uniformly at random then
$\whp$ $\#\Core^*\geq (1-e^{-\Theta(d)})n$.
\end{prop}
\begin{prop} Let $b=\#V(B)$, then the set $J$ defined above satisfies $\#J\geq b/4$
\end{prop}
\begin{Proof} Observe that if $F_B$ is minimal then $\#F_B\leq b+1$. This is because in every
cycle the number of variables equals the number of clauses, and in
the worst case, the path contains at most one more clause than the
number of variables, and the same goes for the chord-cycle. Now
suppose in contradiction that $\#J<b/4$, then there are more than
$3b/4$ variables in $V(B)$, each appearing in at least 6 different
clauses in $F_B$. Thus, $\#F_B>(6\cdot
3b/4)/3=1.5b\underbrace{>}_{b\geq 3}b+1$ (we divided by three as
every clause might have been counted 3 times), contradicting
$\#F_B\leq b+1$.
\end{Proof}\\
The following proposition ``defuses" the dependency between the event
that a bi-cycle (chord-cycle) was included in the graph and the fact
that it doesn't intersect the core variables. In the following
proposition we fix an arbitrary $\pi$ and $\alpha$ in the definition
of $\Core^*$, therefore the probability is taken only over the
randomness in the choice of $\F$.
\begin{prop}\label{CycleProbProp}$Pr[F_{B} \subseteq \F
\text{ and } V(B) \cap \Core = \emptyset]\leq Pr[F_{B} \subseteq
\F]\cdot Pr[J \cap \Core^* = \emptyset].$
\end{prop}
\noindent Before proving Proposition \ref{CycleProbProp}, we establish the following fact.
\begin{lem}\label{BasicIncLemma} For every bi-cycle (chord-cycle) $B$ and
every minimal inducing set $F_{B}$,
$\Core^*(\F,\varphi,\alpha,\pi) \subseteq \Core(\F\cup
F_{B},\varphi,\alpha,\pi)$.
\end{lem}
\begin{Proof}
The lemma is proved using induction on $i$ ($i$ being the
iteration counter in the construction of $\Core$). For the base
case $H'_{0}(\F) \subseteq H_{0}(\F \cup F_B)$, since every variable
in $H'_{0}(\F)$ appears in at most 6 clauses in $F_B$ it holds that
$A_i(\F\cup F_B) \subseteq B_i (\F)$, $i=2,3$. $A_1(\F\cup F_B)
\subseteq B_1(\F)$ holds at any rate as more clauses can only
increase the support, and the set $J$ was not even considered for
$H_0$. Suppose now that $H_{i}'(\F) \subseteq H_{i}(\F \cup
F_{B})$, and prove the lemma holds for iteration $i+1$. If $x \in
H'_{i+1}(\F)$ then $x$ supports at least $d/3$ clauses in which
all variables are in $H'_{i}(\F)$. Since $H'_{i}(\F) \subseteq
H_{i}(\F \cup F_B)$, then $x$ supports at least this number of
clauses with only variables of $H_i(\F \cup F_{B})$. Also, $x$
appears in at most $d/30-6$ clauses with some variable outside of
$H'_{i}(\F)$, again since $H'_{i}(\F) \subseteq H_{i}(\F \cup
F_B)$ and $F_{B}$ contains at most 6 clauses containing $x$, $x$
will appear in no more than $d/30$ clauses each containing some
variable not in $H_{i}(\F \cup F_B)$. We conclude then that $x \in
H_{i}(\F \cup F_B)$.
\end{Proof}\\
This lemma clarifies the motivation for defining $\Core^*$. It is
not necessarily true that $\Core(\F)\subseteq \Core(\F \cup F_{B})$.
For example, a variable which appears in $\Core(\F)$ could
disappear from $\Core(\F \cup F_{B})$ since the clauses in $F_{B}$
make it unstable. Loosely speaking, $\Core^*$ is cushioned enough
to prevent such a thing from happening.\\\\

\begin{Proof}(Proposition \ref{CycleProbProp})
$$Pr[F_B \subseteq \F \text{ and } V(B) \cap \Core = \emptyset] \leq Pr[F_B \subseteq \F \text{ and } J \cap \Core = \emptyset]
= Pr[J \cap \Core = \emptyset | F_B \subseteq \F] Pr[F_B \subseteq
\F].$$ Therefore, it suffices to prove
$$Pr[J \cap \Core = \emptyset | F_B \subseteq \F] \leq Pr[J \cap \Core^* = \emptyset].$$
$$Pr[J \cap \Core^* = \emptyset] = \sum_{F:J \cap \Core^*(F) = \emptyset}Pr[\F=F]\underbrace{\geq }_{Lemma~\ref{BasicIncLemma}}\sum_{F:J \cap \Core(F\cup F_B) = \emptyset}Pr[\F=F]$$ Break each set of clauses $F$ into $F'=F
\setminus F_B$ and $F''=F \cap F_B$, and the latter equals
$$\sum_{F':F' \cap F_B=\emptyset,J \cap \Core(F'\cup
F_B)=\emptyset} \sum_{F'':F''\subseteq F_B} Pr[\F\setminus F_B=F'
\text{ and } \F \cap F_B = F'']$$ Since the two sets of clauses,
$\F\setminus F_B$, and $\F \cap F_B$, are disjoint, and clauses
are chosen independently, the last expression equals,
\begin{align*}
&\sum_{F':F' \cap F_B=\emptyset,J \cap \Core(F'\cup
F_B)=\emptyset}\sum_{F'':F''\subseteq F_B} Pr[\F\setminus
F_B=F']Pr[\F \cap F_B =
F'']= \\
& \sum_{F':F' \cap F_B=\emptyset,J \cap \Core(F'\cup
F_B)=\emptyset} Pr[\F\setminus F_B=F']
\underbrace{\sum_{F'':F''\subseteq F_B} Pr[\F \cap F_B = F'']}_{1}
=
\\
& \sum_{F':F' \cap F_B=\emptyset,J \cap \Core(F'\cup
F_B)=\emptyset}Pr[\F\setminus F_B = F']
\end{align*}
Since $(\F \setminus F_B) \cap F_B =\emptyset$, and clauses are
chosen independently, the event $\{F_B\subseteq \F\}$ is
independent of the event $\{\F\setminus F_B = F'\}$. Therefore,
the latter expression can be rewritten as $$\sum_{F':F' \cap
F_B=\emptyset,J \cap \Core(F'\cup F_B)=\emptyset}Pr[\F\setminus
F_B = F'|F_B \subseteq \F]=
 Pr[J \cap \Core = \emptyset|F_B \subseteq \F].$$
\end{Proof}\\

\begin{lem}\label{lem_chordCycleProbCor} Let $B=B_{2i,k,t}$ be a chord-cycle, then
$Pr[V(B) \cap \Core^* = \emptyset | \# \Core^* = (1-\lambda) n]\leq p(i,k)$
where:\begin{enumerate}
    \item $p(i,k)\leq  \lambda^{(i+\frac{k}{2}-1)/4}$ if $B$
    consists of a $2i$-cycle and a variable-variable $k$-path as a chord.
    \item $p(i,k)\leq  \lambda^{(i+\frac{k}{2})/4}$ if $B$
    consists of $2i$-cycle and a clause-clause $k$-path as a chord.
    \item $p(i,k)\leq  \lambda^{(i+\frac{k-1}{2})/4}$ if $B$
    consists of $2i$-cycle and a variable-clause $k$-path as a chord.
\end{enumerate}
\end{lem}

\begin{Proof}
In $(a)$, we have $i+\frac{k}{2}-1$ variables and $i+\frac{k}{2}$
clauses. To bound the event
$\{J \cap \Core^* = \emptyset\}$, given the size of $\Core^*$, observe that $F_B$ is fixed in
the context of this event, and there is no pre-knowledge whether
$F_B$ is included in $\F$ or not. Therefore, $J$ can be treated as
a fixed set of variables, thus the choice of $\Core^*$ is
uniformly distributed over $J$. Recalling that $\#J \geq
(i+\frac{k}{2}-1)/4$, it follows that
\begin{align*}
Pr[J \cap \Core^* = \emptyset | \# \Core^* = (1-\lambda) n] \leq
\frac{\binom{n-\#\Core^*}{\#J}}{\binom{n}{\#J}}=\frac{\binom{\lambda
n}{(i+\frac{k}{2}-1)/4}}{\binom{n}{(i+\frac{k}{2}-1)/4}} \leq
\lambda^{(i+\frac{k}{2}-1)/4}.
\end{align*}
The last inequality follows from standard bounds on the binomial
coefficients.  This proves $(a)$. In the same manner items $b,c$ are proven (just
counting how many variables and clauses $B$ contains, depending on
the type of its path).
\end{Proof}

\begin{lem}\label{lem_BiCycleProbCor} Let $B=B_{2i,2j,k,t}$ be a bi-cycle, then
$Pr[V(B) \cap \Core^* = \emptyset | \# \Core^* = (1-\lambda) n]\leq p(i,j,k)$
where:\begin{enumerate}
    \item $p(i,j,k)\leq  \lambda^{(i+j+\frac{k}{2}-1)/4}$ if $B$
    consists of a $2i$,$2j$-cycles and a variable-variable $k$-path.
    \item $p(i,j,k)\leq \lambda^{(i+j+\frac{k}{2})/4}$ if $B$
    consists of $2i$,$2j$-cycles and a clause-clause $k$-path.
    \item $p(i,j,k)\leq  \lambda^{(i+j+\frac{k-1}{2})/4}$ if $B$
    consists of $2i$,$2j$-cycles and a variable-clause $k$-path.
\end{enumerate}
\end{lem}
\noindent Lemma \ref{lem_BiCycleProbCor} is proven in a similar way
to Lemma \ref{lem_chordCycleProbCor}.

To complete the proof of Proposition \ref{prop_NonCoreFactorGraph}, we use
the union bound over all possible bi/chord-cycles. We present the proof for the
bi-cycle case with a variable-variable path; the proof for all other cases is identical.
$s=s_{i,j,k}=i+j+\frac{k}{2}-1$ (namely, $\#V(B)=s$ and
$\#F_B=s+1$). We also let $p_\lambda=Pr[\#\Core^*=(1-\lambda) n]$. The probability of $B$ is then at most
\begin{align*}
&\sum_{\lambda=0}^{n} p_\lambda\cdot \sum_{i,j,\frac{k}{2}=1}^{n}\binom{n}{s_{i,j,k}}\cdot
(s_{i,j,k})! \cdot (7n)^{s_{i,j,k}+1}\cdot
\left(\frac{d}{n^2}\right)^{s_{i,j,k}+1} \cdot
\lambda^{s_{i,j,k}/4} \leq\\& \sum_{\lambda=0}^{n} p_\lambda\cdot\sum_{i,j,\frac{k}{2}=1}^{n} 7d\cdot
\left(\frac{7en}{s}\right)^{s} \cdot s^{s} \cdot n^{s+1} \cdot
\left(\frac{d}{n^2}\right)^{s+1} \cdot \lambda^{s/4}  \leq
\sum_{\lambda=0}^{n} p_\lambda\cdot\sum_{i,j,\frac{k}{2}=1}^{n}(7e \cdot d \cdot
\lambda^{1/4})^{s}\cdot \frac{7d}{n} .
\end{align*}
Let us now break the first sum according to the different values of $\lambda$. Proposition \ref{SizeOfHBarPr} implies
that $p_\lambda < n^{-5}$ for all $\lambda > d^{-8}$. Therefore the contribution of this part to the (double) summation
is $O(1/n)$. For $\lambda < d^{-8}$, the latter simplifies to

\begin{equation}\label{eq:lambda}\sum_{\lambda=0}^{n/d^8}p_\lambda\sum_{i,j,\frac{k}{2}=1}^{n}\left(\frac{1}{2}\right)^{s}\cdot
\frac{7d}{n}.
\end{equation}
Finally observe that
$$\sum_{i,j,\frac{k}{2}=1}^{n}\left(\frac{1}{2}\right)^{s}\cdot
\frac{7d}{n} \leq \sum_{i+j+\frac{k}{2} \leq 4\log
n}\frac{7d}{n}+\sum_{i+j+\frac{k}{2}\geq 4\log n}
\left(\frac{1}{2}\right)^{s}\leq (4\log n)^3\cdot
\frac{7d}{n}+n^3\cdot \frac{1}{n^4}=o(1).
$$
Therefore, (\ref{eq:lambda}) sums up to $o(1)$.

\subsection{Proof of Proposition \ref{prop_NonCoreFactorGraph} $(c)$}
The proof is basically the same as that of Proposition
\ref{prop_NonCoreFactorGraph} $(b)$. One defines the same notion of ``cushioned"
core $\Core^*$, and proceeds similarly. We therefore reprove only
the last part -- the union bound over all possible cycles.

First let us bound the number of cycles of length $k$. There are
$\binom{n}{k}$ ways to choose the variables inducing the cycle, and
$k!/2$ ways to order them on the cycle. As for the set of clauses
that induces the cycle, once the cycle is fixed, we have at most
$(7n)^k$ ways of choosing the third variable and setting the
polarity in every clause. We also let $p_\lambda=Pr[\#\Core^*=(1-\lambda) n]$.

Using the union bound, the probability of a cycle of length at least
$k$ in the non-core factor graph is at most
\begin{align*}
&\sum_{\lambda=0}^{n} p_\lambda\cdot\sum_{t=k}^{n}\binom{n}{t}\cdot t! \cdot
(7n)^t\cdot
\left(\frac{d}{n^2}\right)^t \cdot
\lambda^{t/2}   \leq \sum_{\lambda=0}^{n} p_\lambda\cdot\sum_{t=k}^{n}
\left(\frac{7en}{t}\right)^t \cdot t^t \cdot n^t \cdot
\left(\frac{d}{n^2}\right)^t \cdot \lambda^{t/2} \\& =
\sum_{\lambda=0}^{n}p_\lambda \cdot\sum_{t=k}^{ n}(7e \cdot d \cdot \sqrt{\lambda})^t.
\end{align*}
Let us now break the first sum according to the different values of $\lambda$. Proposition \ref{SizeOfHBarPr} implies
that there exists a constant $c$ s.t. $p_\lambda < n^{-3}$ for $\lambda > e^{-cd}$. Therefore the contribution of this part to the (double) summation
is $O(1/n)$. For $\lambda \leq e^{-cd}$, $7e \cdot d \cdot
\sqrt{e^{-cd}}=e^{-\Theta(d)}$.
In this case, the last summation is simply the sum of a
decreasing geometric series with quotient $e^{-\Theta(d)}$, which
sums up to at most twice the first term, namely
$e^{-\Theta(dk)}$. The proposition then follows.

\subsection*{Acknowledgement}

We thank Eran Ofek for many useful discussions and two anonymous referees for many helpful comments. This work was done
while the authors were visiting Microsoft Research, Redmond,
Washington.





\end{document}